\documentclass[aip,rsi,amsmath,amssymb,preprint,nofootinbib]{revtex4-1}
\usepackage{float}
\usepackage{xcolor}
\usepackage{amssymb,amsmath,amsthm}
\usepackage{graphicx}
\usepackage{epstopdf}
\usepackage{subfigure}

\def\d{\mathrm{d}}

\def\eps{{\varepsilon}}
\def\R{\mathbb{R}}
\def\Z{\mathbb{Z}}

\renewcommand{\vec}[1]{\mbox{\boldmath$#1$}}
\def\sech{{\, \mathrm{sech} \, }}

\def\stretch#1{\renewcommand{\arraystretch}{#1}}

\newtheorem{theorem}{Theorem}

\begin{document}

\title{Quantifying the role of folding in nonautonomous flows: the unsteady Double-Gyre}

\author{K.G.D. Sulalitha Priyankara}
\affiliation{Department of Mathematics, Clarkson University, 8 Clarkson Ave, Potsdam, New York, 13699-5815, USA.}
\author{Sanjeeva Balasuriya}
\email{sanjeeva.balasuriya@gmail.com}
\affiliation{School of Mathematical Sciences, University of Adelaide, Adelaide, SA, 5005, Australia.}

\author{Erik Bollt}
\email{bolltem@clarkson.edu}
\affiliation{Department of Mathematics, Clarkson University, 8 Clarkson Ave, Potsdam, New York, 13699-5815, USA.}

\begin{abstract}
We analyze chaos in the well-known nonautonomous Double-Gyre system.   A key focus is on
{\em folding}, which is possibly the less-studied aspect of the 
``stretching + folding = chaos'' mantra of chaotic dynamics.  Despite the Double-Gyre not having the classical
homoclinic structure for the usage of the Smale-Birkhoff theorem to establish chaos, we use the concept of folding to prove the existence of an
embedded horseshoe-map.  We also show how curvature of manifolds can be used to identify
fold points in the Double-Gyre.  This method is applicable to general nonautonomous flows in two dimensions,
defined for either finite or infinite times.
\end{abstract}
\maketitle

\noindent {\bf The well studied Double-Gyre system may be considered as a ``standardized" problem to contrast chaos from mixing. It is often invoked as a model for a chaotic system, and used as a testbed in numerical simulations.  Strangely, however, there does not appear to be a proof in the literature that the system is actually chaotic.  While it is easy to establish that certain stable and unstable manifolds intersect, there are technical impediments in taking the next step to claim that this results in chaos.    By using a new technique we call ``fold re-entrenchment,"  we are able here to show the presence of a Smale horse-shoe embedded in the phase space.  In this process, we  are focussing on the \emph{folding} aspect of chaos, in contrast to highly popular methods such as Finite-Time-Lyapunov Exponents (FTLEs) which target the quantification of  \emph{stretching}.  We further develop a quantification of folding as a system's propensity to develop curvature, and show how this criterion can be highly informative in analyzing the chaotic nature of general systems.   }

\section{Introduction}

A well-known mechanism through which chaos can arise in deterministic dynamical systems is by the combined effect of
stretching and folding.  Stretching will separate nearby points, while folding can abrupt brings together points which were initially
far away.   Various ways which quantify the stretching (most notably finite-time Lyapunov exponents) abound in the 
literature \cite{shadden,Froyland20121612,siam_book,tallapragadaross,117.164502,maouellettebollt}.  {\em Folding}, however, is much less addressed.  In this paper, we specifically focus on the concept of
folding in two different ways.  Firstly, it is used to prove that a highly-studied testbed for numerical methods---the Double-Gyre \cite{shadden}---is chaotic.  While the fact that this system is chaotic is `known' 
anecdotally, it appears that a proof of this fact is not available, and we are able to provide it in this paper using the concept
of folding in a specific way.  Secondly, we propose a method for quantifying folding in general two-dimensional
nonautonomous dynamical systems.  This is through computing the curvature along distinguished one-dimensional curves
of the system.

The presence of stretching and folding in a dynamical system leads to a range of properties usually associated with chaos (sensitivity to initial conditions, presence of countably
many periodic orbits and uncountably many aperiodic ones, the presence of a dense orbit, etc).  Smale's horseshoe map
\cite{ROSSLER1977392,HOLMES19867,ARROWSMITH1993292,Chen2006377,alligoodsaueryorke} forms a paradigm for this mechanism,  and in proving that this system is chaotic the basic strategy
is to exploit the conjugacy of the map's action with shift dynamics on bi-infinite sequences \cite{alligoodsaueryorke,0951-7715-17-4-009}.  Thus,
in proving that two-dimensional maps are chaotic, it is sufficient to establish the existence of horseshoe-like maps
within them.   One standard way in which this arises is through the presence of a transverse intersection between the stable
and the unstable manifold of a fixed point of the map; the Smale-Birkhoff theorem\cite{guckenheimerholmes,holmes,alligoodsaueryorke} provides a method for constructing the horseshoe map in that situation.
The original theorem is for {\em homoclinic} situations; that is, there must be a transverse intersection between the stable 
and the unstable manifold of the {\em same} fixed point of a discrete dynamical system.  The basic intuition is that it is then possible to identify a 
quadrilateral piece of space near the intersection (call it $ A $), which eventually gets mapped back on top of itself exactly 
like a horseshoe map.  The homoclinic
nature is crucial in this argument, since it enables $ A $ to get mapped `all the way round' since after it gets 
pulled out along the unstable manifold direction, it will then get pulled in along the stable manifold direction.

The Smale-Birkhoff theorem does not apply to the Double-Gyre flow \cite{shadden}, since it does not have a
homoclinic structure.   The Double-Gyre was initially proposed by Shadden et.~al.\ \cite{shadden} as a toy model for a
two adjacent oceanic gyres.  It has since taken on a prominent role as a testbed in the development of a range of
numerical diagnostics associated with transport and transport barriers \cite[e.g.]{allshousepeacock,williams,pratt,siam_book,sudharsan,rosi,garaboa,mabollt,mcilhanywiggins,mosovskymeiss,bruntonrowley,lipinskimohseni,ducsiegmund,tallapragadaross, bollt2012measurable, bollt2000controlling, bollt2002manifold, froyland2013}.  Numerics amply demonstrate that there is chaotic transport between the two gyres, which can each
be thought of as a Lagrangian coherent structure \cite{doi:annurev-fluid,PhysRevE.88.013017,2000PhyD..147..352H,Karrasch2015b,Onu201526}.  The field of Lagrangian coherent structures continues to attract
tremendous interest, and there is in particular a multitude of diagnostic techniques that are either being refined or newly
developed for the analysis of fluid transport associated with them.  Well-established methods include finite-time Lyapunov
exponent fields \cite{Shadden2005271,he2016,JGRC:JGRC21415,Nelson201565,johnauc,BozorgMagham2015964,branickiwigi,2016APS..DFD.A8002B}, transfer (Perron-Frobenius) operator approachs \cite{Froyland20091507,PhysRevLett.98.224503,Dellnitz00setoriented,froysanti}, averages along trajectories \cite{wigginsmancho12,Mancho20133530,poje,mancho,mezic} and
curves of extremal attraction/repulsion \cite{blazevski,teramoto,farazmand}.  Other methods include clustering approaches \cite{PhysRevE.93.063107,JGRC:JGRC21415,froylandgelhe1},
topological entropy \cite{balibrea,PMID:27036190,TUMASZ2013117}, ergodic-theory related approaches \cite{Budišić20121255} and curvature \cite{mabolltcurvature,maouellettebollt}.  The latter approach is
particularly relevant to the current paper, and will be revisted later.  The main point, though, is that the Double-Gyre is often
used to {\em test} these methods, sometimes against each other \cite{allshousepeacock}.  In doing so, the `complicated'
(i.e., chaotic) nature of the Double-Gyre is taken as given.  However, there is as yet no proof that it is actually chaotic!

From the theoretical perspective, the impediment to using the Smale-Birkhoff theorem is that the entity separating the
two gyres is not homoclinic, but rather {\em heteroclinic}.  That is, it is associated with the stable manifold of a fixed point
(of a relevant Poincar\'e map), and the unstable manifold of a {\em different} fixed point, intersecting.  The standard horseshoe
construction fails in this situation.  An approach might be to appeal to an extension of the Smale-Birkhoff theorem due to
Bertozzi \cite{bertozzi}, in which she considers a `heteroclinic cycle' in which intersection patterns between 
stable/unstable manifold structures of a collection of fixed points forms a cycle.  Under {\em generic} conditions, it is then
shown that a horseshoe construction can be made in this situation as well \cite{bertozzi}; effectively, the region $ A $ gets
mapped around, going near each fixed point, and eventually returning to form a horseshoe-like set falling on top of $ A $.
Unfortunately, the Double-Gyre does {\em not} fall into this generic situation.  While there is a heteroclinic cycle geometry
in the Double-Gyre, only one of the connections between fixed points possesses the generic transverse intersection pattern.
All other connections are situations in which a stable manifold {\em coincides} with an unstable manifold, and thus the
heteroclinic extension \cite{bertozzi} to the Smale-Birkhoff theorem is inapplicable.  

Given the importance of the Double-Gyre as a testbed, and the implicit agreement that it {\em is} chaotic, an actual
{\em proof} of its chaotic nature would seem important.  We provide exactly that in Section~\ref{sec:horseshoe}.  We first
develop analytical approximations to the stable and unstable manifolds.  These are then used to identify `fold points' which
are the basis for a horseshoe construction, leading to Theorem~\ref{modeltent} in which we establish the chaotic nature
of the Double-Gyre.

A main ingredient leading to chaos appearing in the Double-Gyre is the fact that the stable and unstable manifolds {\em fold}.
We address this issue in a complementary fashion in Section~\ref{sec:curvature}.  Here, we are inspired by recent work on using
{\em curvature} in Lagrangian coherent stucture analysis \cite{mabolltdiff}.   In the current context, though, the argument is simple:
if stable/unstable manifolds fold, then the curvature at those fold points must get anomalously large. Using the Double-Gyre
as a testbed, we both numerically and theoretically track such points of large curvature.  We establish numerically that the
fold points do indeed possess the behavior established in our proof of chaos in the Double-Gyre.  Using the curvature in this
way can be done for general two-dimensional nonautonomous flows.  The Double-Gyre is time-periodic, which allows for
thinking of the dynamical system either in continuous time, or in discrete time (in relation to a Poincar\'e map).  However, 
it is possible to use the curvature in nonautonomous systems with any time-dependence, by thinking of the stable and
unstable manifolds as being attached to hyperbolic trajectories \cite{jusmallwiggins,850647,esvan} rather than fixed points.  Moreover, using curvature
in this way can also be done for specialized curves arising from using any diagnostic procedure in {\em finite-time} flows.

\section{Horseshoe map chaos in the Double-Gyre flow}
\label{sec:horseshoe}

The Double-Gyre flow was initially introduced by Shadden et.~al. \cite{shadden}, and has since been studied extensively as a canonical example
of complicated transport in nonautonomous flows 
\cite[e.g.]{allshousepeacock,williams,pratt,siam_book,sudharsan,rosi,garaboa,mabollt,mcilhanywiggins,mosovskymeiss,bruntonrowley,lipinskimohseni,ducsiegmund}.
Its flow is given by
\begin{equation}
\stretch{1.8} \left. \begin{array}{ll}
\dot{x}_1 &= - \pi A \sin \left[ \pi \phi(x_1,t) \right] \cos \left[ \pi x_2 \right] \\
\dot{x}_2 &= \pi A \cos \left[ \pi \phi(x_1,t) \right] \sin \left[ \pi x_2 \right] \frac{\partial \phi}{\partial x_1}(x_1,t)
\end{array} \right\} \, , 
\label{eq:doublegyre}
\end{equation}
in which $ A > 0 $ and $ 0 < \eps \ll 1 $, and
\[
\phi(x_1,t) := \eps \sin \left( \omega t \right)  x_1^2 + \left( 1 - 2 \eps \sin \left( \omega t \right) \right) x_1 \, .
\]
This is usually viewed in the spatial domain $ \Omega := [0,2] \times [0,1] $, and when $ \eps = 0 $ possesses two
counter-rotating gyres: one in $ (0,1) \times (0,1) $ and the other in $ (1,2) \times (0,1) $.  This is a steady
situation in which the gyres are separated by a heteroclinic manifold $ x_1 = 1 $, which is the stable manifold of $ (1,0) $
and the unstable manifold of $ (1,1) $.  This manifold can be expressed parametrically by
\begin{equation}
\bar{x}_1(t) = 1 \quad , \quad \bar{x}_2(t) = \frac{2}{\pi} \, \cot^{-1} e^{\pi^2 A t} \, \quad ; \quad t \in \R \, ,
\label{eq:heteroclinic}
\end{equation}
which is an exact solution to (\ref{eq:doublegyre}), with $ t $ representing time, when $ \eps = 0 $.  Here, $ t = 0 $ 
corresponds to $ x_2 = 1/2 $, and $ \bar{x}_2(t) \rightarrow 1 $ when $ t \rightarrow - \infty $ and $ \bar{x}_2(t) \rightarrow 0 $
when $ t \rightarrow \infty $.  

When $ \eps \ne 0 $, the flow (\ref{eq:doublegyre}) is nonautonomous.  In this case of the `classical' Double-Gyre,
it is time-periodic as well (for an analysis similar to what is to be presented for the {\em aperiodic} Double-Gyre, see
\cite{siam_book}).  Despite the nonautonomous nature, the lines $ x_1 = 0 $, $ x_1 = 2 $, $ x_2 = 0 $ 
and $ x_2 = 1 $ (i.e., the boundary of $ \Omega $) remain invariant.  Thus,
there is no possibility of chaotic motion being created in the system by the mechanism reported by
Bertozzi \cite{bertozzi}, which requires the heteroclinic network to break apart all the way around.  
However---as is well-known anecdotally and numerically though
a proof does not seem to appear in the literature yet---the heteroclinic manifold which
separates the two gyres, {\em does} break apart such that transverse intersections are created.  The proof is straightforward.

\begin{theorem}[Heteroclinic intersections]
\label{theorem:intersect}
There exists $ \eps_0 $ such that for $ \left| \eps \right| \in (0,\eps_0) $, at any time $ t $, the stable and unstable manifolds adjacent to $ x_1 = 1 $ intersect 
each other transversely an infinite number of times.
\end{theorem}

\begin{proof}
See Appendix~\ref{sec:intersect}.
\end{proof}

Despite perhaps conventional belief, Theorem~\ref{theorem:intersect} does not in and of itself prove the presence of
{\em chaos} in the Double-Gyre.  The original Smale-Birkhoff Theorem \cite{guckenheimerholmes,holmes,alligoodsaueryorke} can only guarantee
chaos, in the sense of symbolic dynamics, for {\em homo}clinic tangles.  If the heteroclinics on the 
outer boundaries of $ \Omega $
also exhibited heteroclinic tangling, then the results of Bertozzi \cite{bertozzi} can help extend this result.  This
is because fluid would travel from one heteroclinic tangle to the next, and so on, until arriving back again in the region of the
initial tangle.  It can be shown \cite{bertozzi} that dynamics similar to Smale's horseshoe map \cite{guckenheimerholmes,holmes,alligoodsaueryorke} ensues, and a continual repetition of this process can be proven to produce chaotic dynamics.  However, in this case
the heteroclinic manifolds on the boundary of $ \Omega $ in the Double-Gyre {\em do not break apart}.  Indeed, (\ref{eq:doublegyre}) was proposed \cite{shadden} to preserve
these boundaries, in order for it to be a model for an oceanic Double-Gyre enclosed by land.  Therefore, additional work is needed to establish
how {\em chaotic} transport occurs in the Double-Gyre due to the heteroclinic tangle adjacent to $ x_1 = 1 $.

The crux to the argument is the fact that the stable and unstable manifolds in the heteroclinic tangle have {\em folds} in them.
We will show that fluid adjacent to such folds gets transported around the gyres and back again into the heteroclinic tangle region.
In doing so, we will need analytical approximations for the stable and unstable manifolds, and the hyperbolic trajectories
to which they are attached, for small $ \left| \eps \right| $.  In the following, we think of these entities as {\em nonautonomous}
ones, i.e., not necessarily in terms of a Poincar\'e map.  From this viewpoint, a hyperbolic trajectory is defined in terms
of exponential dichotomy conditions \cite{coppel,battellilazzari,palmer,siam_book,eigenvector}, and its local stable manifold is associated
with the projection operator of the exponential dichotomy.  The global stable manifold is of course the continuation of this.
All these entities are therefore parametrized by time $ t  \in \R $.  The time-periodicity property of the Double-Gyre
will allow for identification of these nonautonomous entities equivalently in terms of a Poincar\' map $ P_t $ which takes the flow from
time $ t $ to $ t + 2 \pi / \omega $; the hyperbolic trajectory location would be a hyperbolic fixed point of $ P_t $ and the nonautonomous stable manifold will coincide with the stable manifold (with respect to $ P_t $) of this hyperbolic fixed point.  The
advantage of the nonautomous viewpoint is that the variation with $ t $ is retained, whereas if considering a Poincar\'e map
$ P_t $ then it is necessary to think of $ t \in [0, 2 \pi/\omega) $.  Thus, there is in actuality a {\em family} of Poincar\'e maps.  
We will go back and forth between these continuous-time and discrete-time viewpoints, as needed.

Using the continuous-time approach, the hyperbolic trajectory and (a part of) its stable manifold can be approximated by the following
theorem. 

\begin{theorem}[Stable manifold]
\label{theorem:stable}
Let $ P \in \R $ be given.  Then, there exists $ \eps_0 $ such that for $ \left| \eps \right| \in (0,\eps_0) $, the saddle point at $ (1,0) $ when $ \eps = 0 $, perturbs to a time-varying hyperbolic trajectory $ \vec{x}_h^s(t) = \left( x_h^s(t), 0 \right) $  given by
\begin{equation}
x_h^s(t) = 1 + \eps \cos \theta \sin \left( \omega t + \theta \right)  + {\mathcal O}(\eps^2) \quad ; \quad 
\theta := \tan^{-1} \frac{\omega}{A \pi^2} \, .
\label{eq:xshyp}
\end{equation}
Moreover, the stable manifold emanating from  $ \vec{x}_h^s(t) $ at a time $ t $ can be approximated in the vicinity of $ x_1 = 1 $ in the parametric form
\begin{small}
\begin{equation}
\stretch{2} \left. \begin{array}{ll}
\displaystyle x_1^s(p,t) & \displaystyle \! \! \! \! = 1 \! + \! \eps \frac{ \pi^2 A}{\sech \left( \pi^2 A p \right) }
\! \int_p^\infty \!  \! \! \tanh \left( \pi^2 A \tau \right)
\! \sech \left( \pi^2 A \tau \right)  \sin \left[ \omega (\tau \! + \! t \! - \! p) \right]
  \d \tau  \! + \! {\mathcal O}(\eps^2) \\
 \displaystyle x_2^s(p,t) & \displaystyle = \frac{2}{\pi} \, \cot^{-1} e^{\pi^2 A p}
 \end{array} \right\}  \, , 
\label{eq:stable}
\end{equation}
\end{small}
for $ p \in [P,\infty) $, and moreover if its reciprocal slope at $ \vec{x}_h^s(t) $, is $ \theta_s(t) $, then
there exists $ K_s $ such that $ \left| \theta_s(t) \right| \le \eps^2 K_s $ for $ (t,\eps) \in \R \times [0,\eps_0) $.
\end{theorem}

\begin{proof}
See Appendix~\ref{sec:stable}.
\end{proof}

An alternative expression for the leading-order stable manifold can be obtained by eliminating the parameter $ p $ from
(\ref{eq:stable}).  Since 
\[
p = \frac{1}{\pi^2 A} \ln \left( \cot \frac{\pi x_2}{2} \right) \, 
\]
the stable manifold's leading-order term at each time $ t $ can be expressed as a graph from $ x_2 $ to $ x_1 $.  Now,
a similar theorem holds for the unstable manifold:

\begin{theorem}[Unstable manifold]
\label{theorem:unstable}
Let $ P \in \R $ be given.  Then, there exists $ \eps_0 $ such that for $ \left| \eps \right| \in (0,\eps_0) $, 
the saddle point at $ (1,1) $ when $ \eps = 0 $, perturbs to a time-varying hyperbolic trajectory $ \vec{x}_h^u(t) = \left( x_h^u(t), 1 \right) $, given by
\begin{equation}
x_h^u(t) = 1 + \eps \cos \theta \sin \left( \omega t - \theta \right)  + {\mathcal O}(\eps^2) \quad ; \quad 
\theta := \tan^{-1} \frac{\omega}{A \pi^2} \, .
\label{eq:xuhyp}
\end{equation}
Moreover, the unstable manifold emanating from  $ \vec{x}_h^u(t) $ at a time $ t $ can be approximated in the vicinity of $ x_1 = 1 $ in the parametric form
\begin{small}
\begin{equation}
\stretch{2} \left. \begin{array}{ll}
\displaystyle x_1^u(p,t) & \displaystyle \! \! \! \! = 1 \! - \! \eps \frac{ \pi^2 A}{\sech \left( \pi^2 A p \right) }
\! \int_{-\infty}^p \!  \! \! \tanh \left( \pi^2 A \tau \right)
\! \sech \left( \pi^2 A \tau \right)  \sin \left[ \omega (\tau \! + \! t \! - \! p) \right]
  \d \tau  \! + \! {\mathcal O}(\eps^2) \\
 \displaystyle x_2^u(p,t) & \displaystyle = \frac{2}{\pi} \, \cot^{-1} e^{\pi^2 A p}
 \end{array} \right\}  \, , 
\label{eq:unstable}
\end{equation}
\end{small}
for $ p \in (-\infty, P] $, and moreover if its reciprocal slope at $ \vec{x}_h^u(t) $ is $ \theta_u(t) $, then
there exists $ K_u $ such that $ \left| \theta_u(t) \right| \le \eps^2 K_u $ for $ (t,\eps) \in \R \times [0,\eps_0) $.
\end{theorem}

\begin{proof}
The proof is similar to that of Theorem~\ref{theorem:stable}, and will be skipped.
\end{proof}

\begin{figure}[t]
\includegraphics[width=0.5 \textwidth, height=0.25 \textheight]{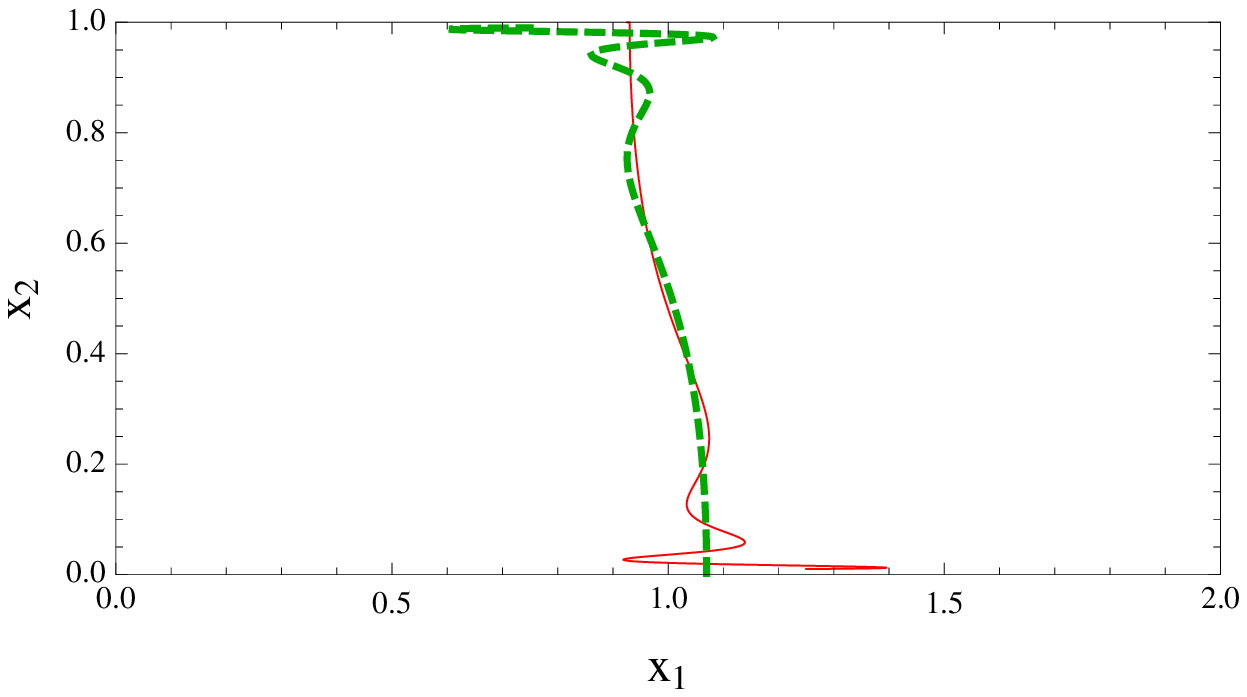}
\includegraphics[width=0.5 \textwidth, height=0.25 \textheight]{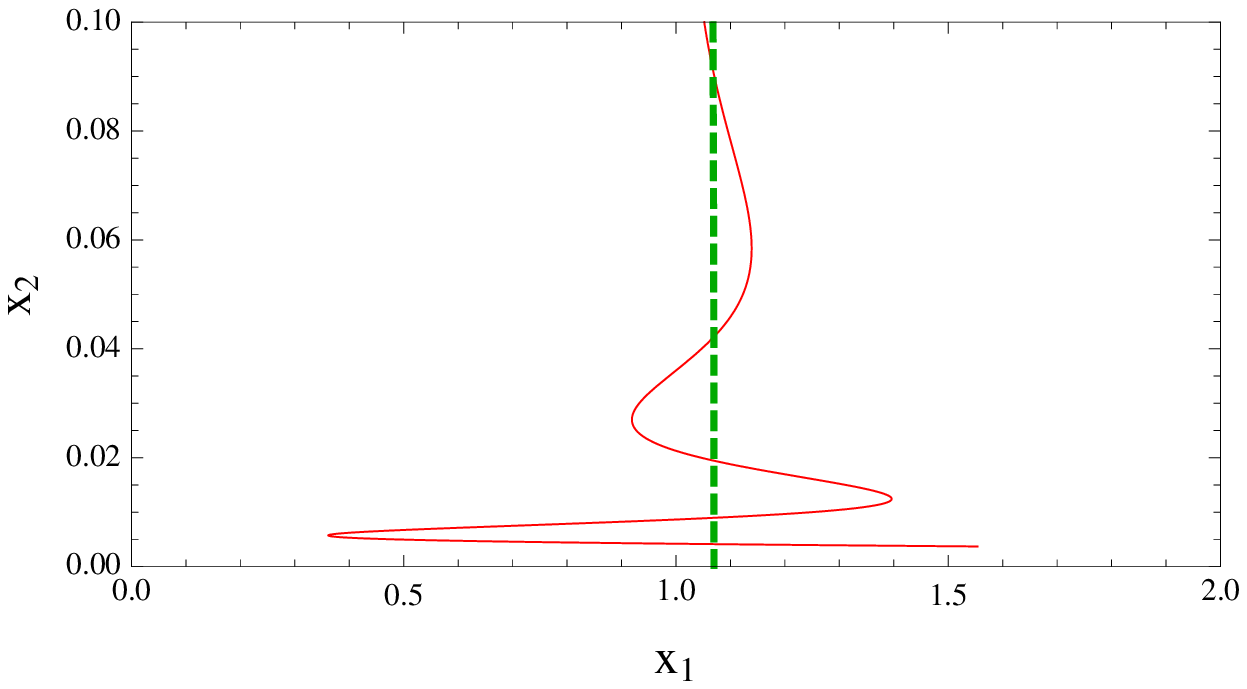}
\caption{Intersection pattern of stable (dashed green) and unstable (solid red) manifolds predicted by Theorems~\ref{theorem:stable} and
\ref{theorem:unstable}, obtained using the $ {\mathcal O}(\eps) $-formul\ae{} at $ t = 0 $ with $ \eps = 0.3 $, $ A = 1 $
and $ \omega = 40 $: in the full domain (left) and zoomed in close to $ x_2 = 0 $ (right).}
\label{fig:manifolds}
\end{figure}

By taking the limit as $ p \rightarrow \infty $ of the $ p $-derivative of the expression~(\ref{eq:stable}), it is possible to
show that the direction of emanation of the stable manifold remains vertical to $ {\mathcal O}(\eps) $.  The same is true
for the unstable manifold; these observations are a special case of the manifold emanation theory developed in
\cite{eigenvector}.  Now, we have already established that the unstable and stable manifolds intersect infinitely often.  Using the expressions in
Theorems~\ref{theorem:stable} and \ref{theorem:unstable}, the nature of this intersection pattern, and the lobes created
as a result of these intersections, can be determined.
We show the intersection pattern at a particular time instance in Fig.~\ref{fig:manifolds}, which was produced with the analytical approximation obtained above, but the computation of the unstable manifold was
stopped after a point.  The unstable manifold can be represented as $ x_1=x_1(x_2) $ for $ p < P_m $ (where $ P_m $ is
an unspecified value), because for $ p \rightarrow
- \infty $, the unstable manifold approaches the hyperbolic trajectory $ \vec{x}_h^u $, from which the unstable manifold
emanates in a well-defined manner.  In this region, we shall refer to the unstable manifold as the {\em primary} unstable manifold, for which (\ref{eq:unstable})
gives a good approximation for small enough $ \left| \eps \right| $.   Larger $ p $-values corresponds to approaching $ x_2 = 0 $,
and here, the unstable manifold will criss-cross the stable manifold infinitely often between the displayed ending and $ x_2 = 0 $.   The stable
manifold near $ x_2 = 0 $ is nearly a straight line emanating upwards from the point $ \vec{x}_h^s(t) $.   However,  the intersection points with the criss-crossing unstable manifold must 
accumulate to $ \vec{x}_h^s(t) $, forcing the corresponding lobes to get elongated in the $ \pm x_1 $-directions in order to maintain incompressibility.
Thus, the unstable manifold in this region will be influenced by global effects, and hence the expression (\ref{eq:unstable}) becomes illegitimate.  It may not be possible to represent the unstable manifold in the form $ x_1 = x_1(x_2) $ in this {\em non-primary} region.
We are able to prove that this is indeed the case, while highlighting a particular behavior.

\begin{figure}[t]
\centering
\includegraphics[width=0.7 \textwidth]{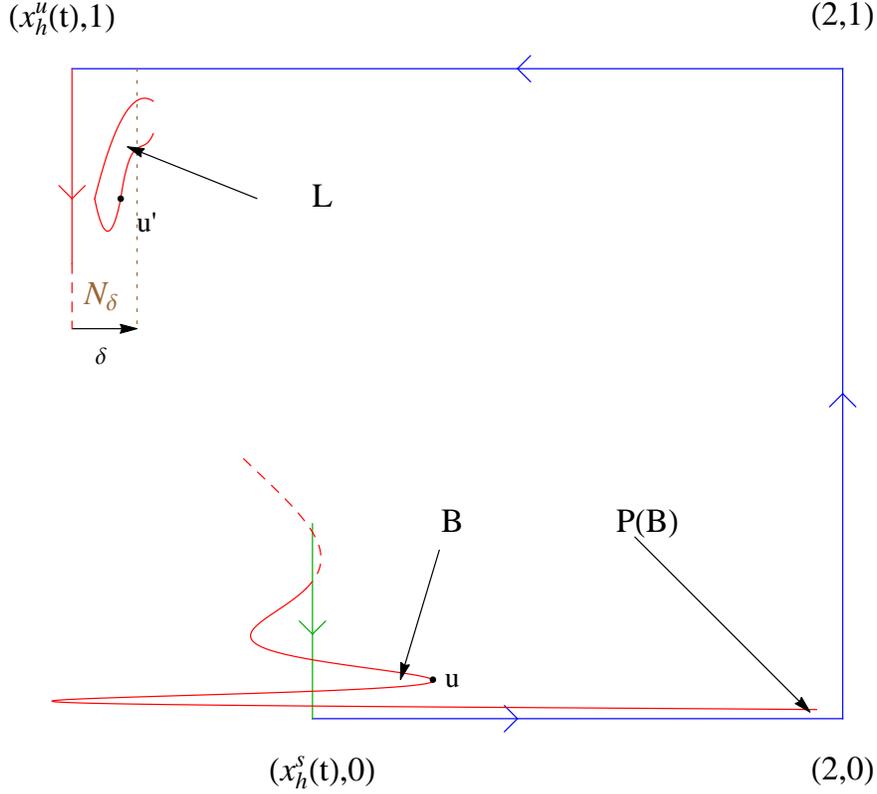}
\caption{The geometry around the right gyre which ensures that the unstable manifold returns to within $ \delta $ of itself after wrapping around
the boundary of $ \Omega $, as described in Theorem~\ref{theorem:reentrench} and Appendix~\ref{sec:reentrench}.}
\label{fig:reentrench}
\end{figure}

\begin{theorem}[Fold re-entrenchment]
\label{theorem:reentrench}
Let $ t \in \R $, and suppose $ \delta > 0 $ is given.  Define  $ N_\delta $ be the one-sided neighborhood of the primary unstable manifold of width
$ \delta $, near the hyperbolic trajectory location $ \vec{x}_h^u(t) = \left( x_h^u(t), 1 \right) $, as shown in Fig.~\ref{fig:reentrench}.
Then, there exists $ \eps_0 $ such that for any $ \eps \in (0,\eps_0) $, the unstable manifold emanating from $ \vec{x}_h^u(t) $
will wrap around the gyre and re-enter $ N_\delta $, forming a fold in the
sense that there is a region within this neighborhood such that a horizontal line will intersect the unstable manifold at least twice.
\end{theorem}

\begin{proof}
See Appendix~\ref{sec:reentrench}.
\end{proof}

The geometry associated with Theorem~\ref{theorem:reentrench} is shown in Fig.~\ref{fig:reentrench}.  There is a
heteroclinic network (shown in blue) connecting the nonautonomous hyperbolic trajectories $ \vec{x}_h^s(t) $ and $
\vec{x}_h^u(t) $ along the outer boundaries of $ \Omega $.  This figure only shows the network around the right
gyre, but there is a similar one around the left.  We note that this is a degenerate situation in that the heteroclinic network 
does {\em not} break apart in
a transverse way.  The parts along the boundary of $ \Omega $ simply persist as straight lines.  This is to be contrasted
with the results 
of Bertozzi \cite{bertozzi}  that {\em generically}, heteroclinic networks which exist for $ \eps = 0 $ break apart through
transverse intersections along each heteroclinic segment.  The Double-Gyre does not follow this, because the nature of
the flow is such that the boundary of $ \Omega $ is forced to remain invariant and regular.  Therefore, Bertozzi's method
for proving existence of a chaotic Smale-horseshoe and chaotic transport in a heteroclinic tangle does not apply for the Double-Gyre.  This is because of this reason that we have had to establish
Theorem~\ref{theorem:reentrench} as a first step in our alternate proof of chaos.  

The main point of Theorem~\ref{theorem:reentrench} is that an unstable manifold segment $ L $ with a fold in it can be found
in any arbitrarily small strip of width $ \delta $ near the primary part of the unstable manifold emanating from $  \left( x_h^u(t),1 
\right) $.  The precise shape of this unstable manifold segment is unknown; for example, it may possess many folds.  However, Theorem~\ref{theorem:reentrench} ensures that there will be {\em at least one} fold, in the sense that on the two sides of 
such a fold point, the unstable manifold
has a larger $ x_2 $ value than at the fold point.  We note that there is no claim that fold points are mapped to fold points.  That is,
it is not necessarily true that the point labeled $ u $ in Figure~\ref{fig:reentrench} will eventually flow to the leading fold in $ L $.  
Its image, $ u' $, need not be a fold point at all.

Now with the re-entrenchment theorem, we are ready to show that the double-gyre flow has an embedded horseshoe map.  The standard Smale horseshoe map is well know \cite{robinson,mybook} to be the map of rectangle, $T_\gamma:R\rightarrow R$ across itself, which in  briefest terms, implies the standard package of results corresponding to fully developed chaos.  Here, our ``rectangle'' will be
slightly different.

\begin{theorem}[Horseshoe Map] 
The Double-Gyre system has an embedded horseshoe, near the point $(1,1)$.  As such, the dynamics of the system is equivalent to a shift-map on a restricted subset, and there is fully developed chaos, at least on this subset.
\label{modeltent}
\end{theorem}

\begin{figure}[t]
\centering
\subfigure[]{\includegraphics[width=0.65 \textwidth]{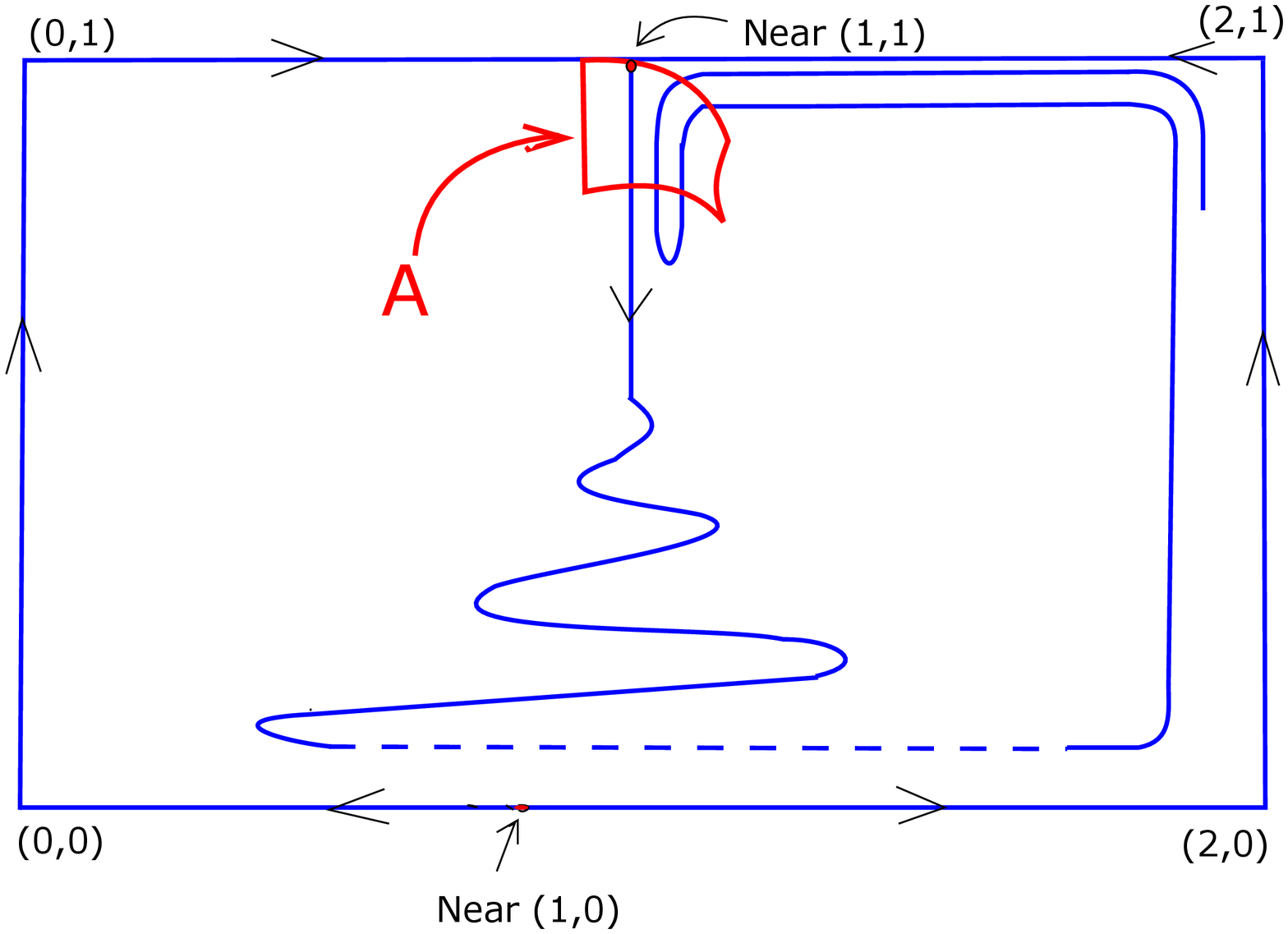}}
\hspace{0.15in}
\subfigure[]{\includegraphics[width=0.3 \textwidth]{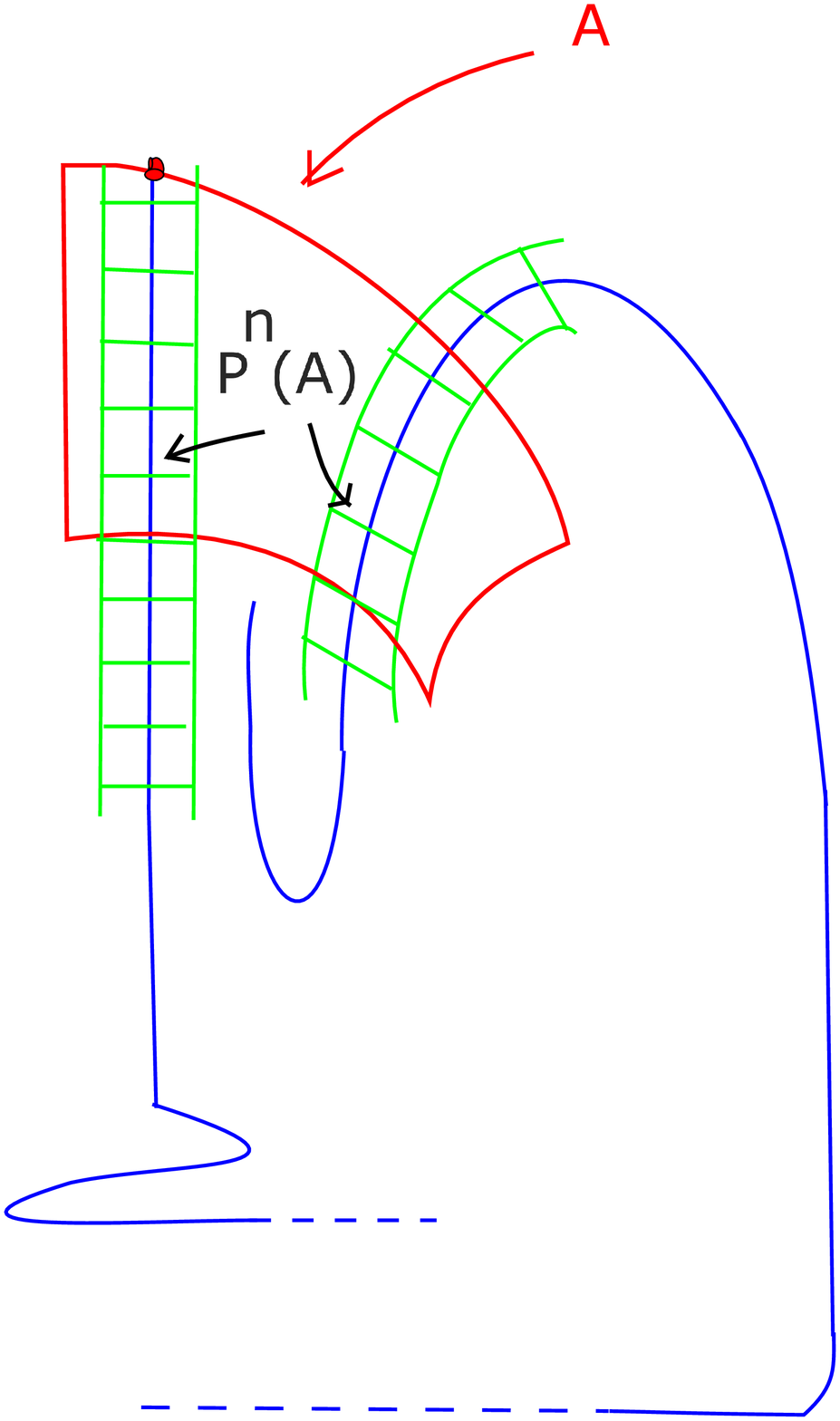}}
\caption{
A topological horseshoe embedded in the Double-Gyre.  (a) As described in proof of Theorem \ref{modeltent}, a topological rectangle set labelled $A$ and shown in red, can be defined transversally to the re-entrenchment region.  (b)  There is a time $n>0$ such that $P^n(A) \cap A$ has stretched into two branches.
}
\label{fig:horsepic}
\end{figure}

\begin{proof}

See Fig.~\ref{fig:horsepic}.  Near the point, $(1,1)$, the unstable manifold shown has been established above.  A set $ A $,
shown by the red boundaries, will be constructed in the re-entrechment region guaranteed by Theorem~\ref{theorem:reentrench}.  A ``vertical'' line, parallel to the emergent unstable manifold, comprises its left boundary.  We next note that there is an infinite number of re-entrenching lobes
accumulating to the primary unstable manifold.  Thus, the curves associated with these lobes, while entering the region
``horizontally,'' will become ``vertical'' in approaching the unstable manifold.  This enables the drawing of the ``top'' boundary of
$ A $ as a curve which
passes through the hyperbolic trajectory but is then normal to each of the curve segments comprising the re-entrenching lobes;
see Fig.~\ref{fig:horsepic}.  The ``bottom'' boundary can also be constructed using the same orthogonality idea.  Finally, the
``right'' boundary of $ A $ is formed by drawing a curve which does not intersect any lobe.  Having constructed $ A $, choose
$ n $ to be large enough such that when applying the strobing Poincar\'e map $P$ $ n $-times to $ A $, part of the set $ A $ will stretch
along the unstable manifold {\em and re-entrench}.   Meanwhile, since $ A $ also ``shrinks'' towards the unstable manifold by the
action of $ P^n $, there will continue to be a part of $ P^n(A) $ which remains within $ A $.  Therefore, the set $ P^n(A) \cap A $
will consist of at least two strips as shown on the right in Fig.~\ref{fig:horsepic}.  While $ A $ is {\em not} a rectangle as in
the usual horseshoe construction \cite{alligoodsaueryorke,robinson,mybook}, this process generates 
an embedded horseshoe \cite{robinson,mybook}.  Define the map $T=P^n$, and let $\Gamma=\cap_{i=-\infty}^\infty T^i(A)$.  Then $T:\Gamma \rightarrow \Gamma$ is semi-conjugate to a Bernoulli shift map on two-symbols, $s:\Sigma_2\rightarrow \Sigma_2$, with details in standard references.  We have claimed only semi-conjugacy since without showing uniform contraction, then it is possible that many points are symbolized by one symbolic sequence.
\end{proof}

Notice that this presentation of an embedded horseshoe is by direct construction, rather than the usual Smale-Birkhoff theorem that follows showing a transverse intersection of stable and unstable manifold, which fails for reasons already described.  Instead we have relied largely on the re-entrenchment theorem.

\section{Folding defined by curvature in the Double-Gyre flow}
\label{sec:curvature}

We have established the existence of chaos in the Double-Gyre system for small enough $ \eps $.  The crux of
this argument comes from the lobes re-entrenching.  Now, these lobes are specifically formed though the {\em folding of
the manifolds}.  The relevance of folding is less studied that stretching (for which, for example, finite-time Lyapunov
exponents \cite{shadden,allshousepeacock,bruntonrowley,lipinskimohseni,garaboa,neighborhood} are a valuable tool), though both contribute toward chaotic transport.  In this section, we follow a recently
emerging idea \cite{maouellettebollt,gajamannage} of examining the folding process in terms of curvature of the manifolds.    Specifically,
we follow the points of high curvature in determining where the ``folding is generated," and the ``stretching'' of the regions
in-between.  Thus, we highlight how stretching and folding interplay in generating the horseshoe-driven chaotic motion
in the Double-Gyre.  We use both
analytical and numerical methods in this analysis, and obtain similar results.

\begin{figure}[t]
\subfigure[$ x_1^s(p) $]{\includegraphics[width=0.5 \textwidth, height=0.25 \textheight]{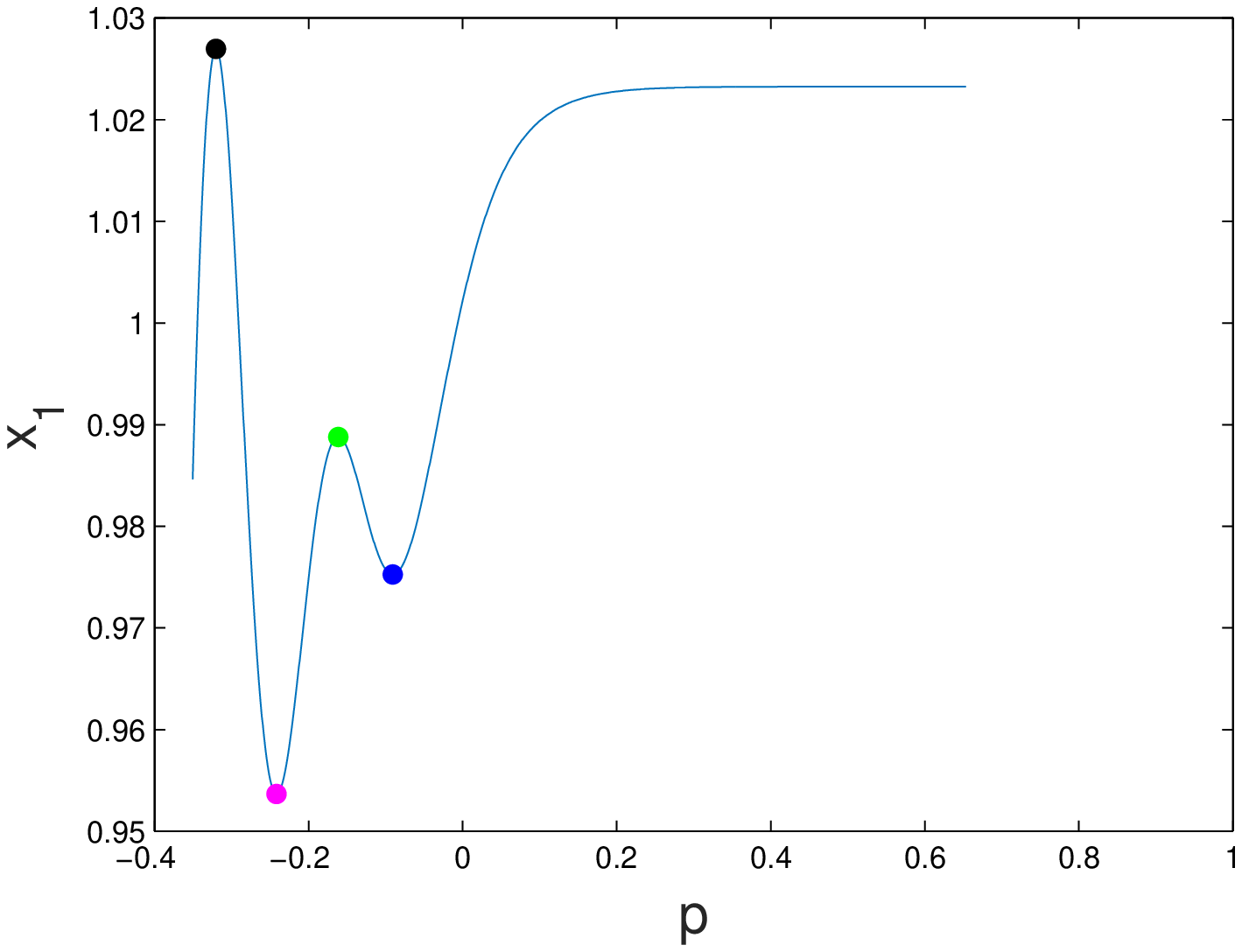}}
\subfigure[$ x_2^s(p) $]{\includegraphics[width=0.5 \textwidth, height=0.25 \textheight]{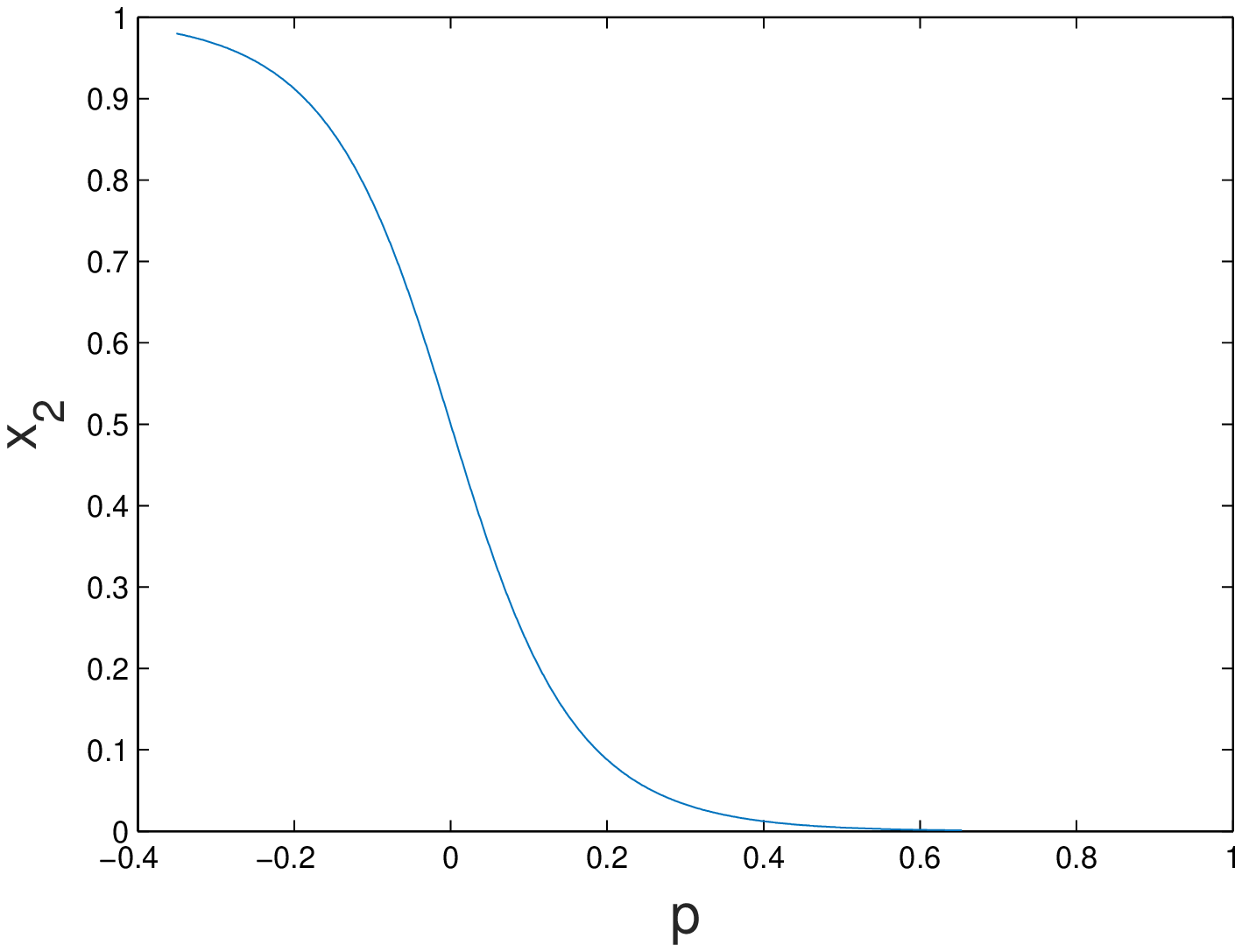}}
\caption{Identification of the first four fold points in the stable manifold, with parameters $A=1, \omega=40, \eps=0.1$ and 
$t=0$.}
\label{fig:foldpointsp}
\end{figure}

\begin{figure}[t]
\subfigure[Full space]{\includegraphics[width=0.5 \textwidth, height=0.25 \textheight]{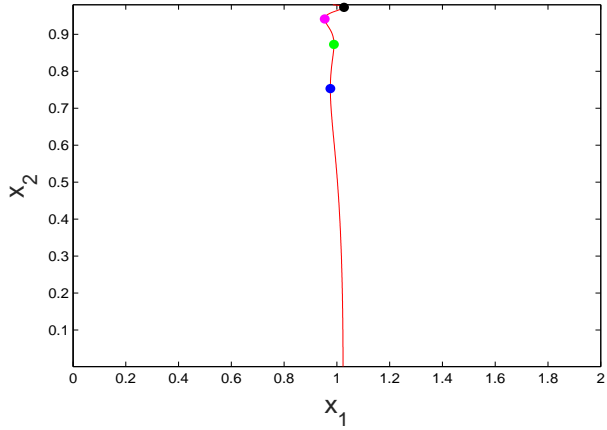}}
\subfigure[Close-up]{\includegraphics[width=0.5 \textwidth, height=0.25 \textheight]{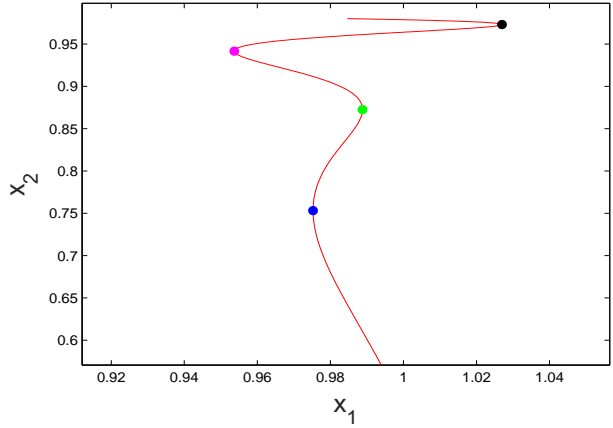}}
 \caption{The four fold points determined in Fig.~\ref{fig:foldpointsp}, illustrated in $ (x_1,x_2) $-space.}
 \label{fig:foldpoints}
 \end{figure}

The analytical expressions in (\ref{eq:xshyp}) and (\ref{eq:xuhyp}) allow for a $ {\mathcal O}(\eps) $ parametric representation 
of the primary segments of the stable and unstable manifolds, in terms of the parameter $ p $, at each fixed time $ t $.
These expressions enable the determination of the location of fold points, distance between points on each manifold, and 
also the curvature at each point on the manifold, as shown in Appendix~\ref{sec:approx}.  Since $ x_2^s(p) $ is monotonic
in $ p $, fold points can simply be obtained by examining turning points of $ x_1^s(p) $ with respect to $ p $; these also 
represent turning points with respect to the variable $ x_2^s $.  We show in Fig.~\ref{fig:foldpointsp} the first four
fold points of the stable manifold \footnote{Bearing in mind that $ p \rightarrow \infty $ approaches the hyperbolic trajectory
$ \vec{x}_h^s(t) $, these correspond to the largest $ p $ values for which $ d x_1^s/d p $ is zero.}, as shown by the
dots in (a).  The locations of these in $ (x_1,x_2) $ space are shown in Fig.~\ref{fig:foldpoints}, where (b) presents a closeup
view of (a).  The same color-coding is used for the four points in both Fig.~\ref{fig:foldpointsp} and \ref{fig:foldpoints}.  We
note that, because the stable manifold curve must intersect the unstable manifold curve (not shown in Fig.~\ref{fig:foldpoints},
but this emanates downwards from near $ (1,1) $) infinitely many times, there must be infinitely many fold points.  We 
only show the first four, since by Theorem~\ref{theorem:stable} the approximation (\ref{eq:stable}) breaks down in the limit $ p \rightarrow - \infty $.  This is because  the stable manifold extends outwards and is impacted by swirling around the boundaries of the Double-Gyre, whereas
the expression (\ref{eq:stable}) is only locally valid near $ x_1 = 1 $.

\begin{figure}[tb]
\centering
\subfigure[]{\includegraphics[width=0.49 \textwidth, height=0.25 \textheight]{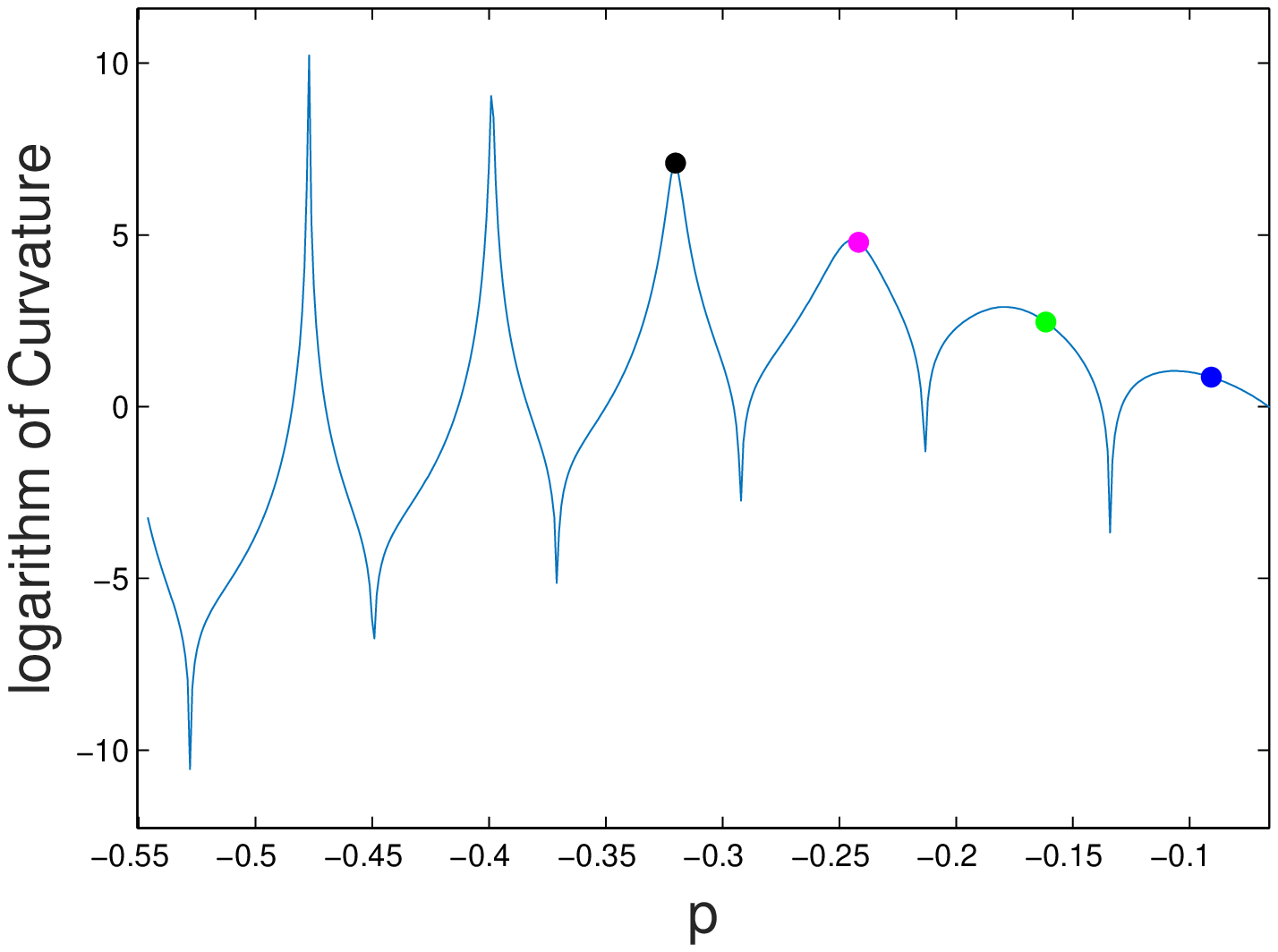}} 
\subfigure[]{\includegraphics[width=0.49 \textwidth, height=0.25 \textheight]{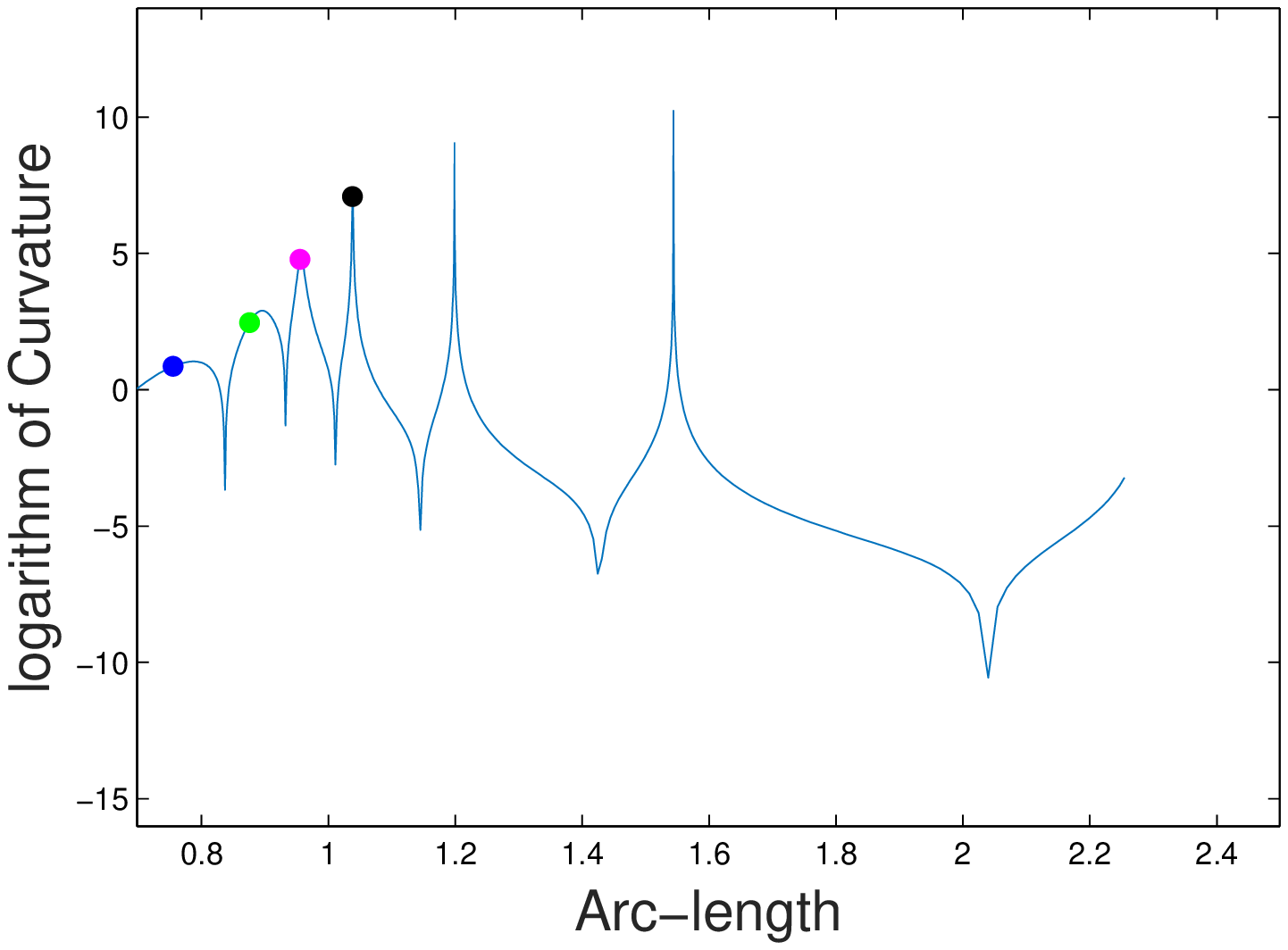}}
 \caption{The logarithm of the curvature along the stable manifold plotted against (a) $ p $, and (b) arclength,
 at the same parameter values as in Figs.~\ref{fig:foldpointsp} and \ref{fig:foldpoints}.}
\label{fig:curvature}
 \end{figure}

In this case we have worked with analytical expressions, and have the advantage of knowing that the turning points of
$ x_1^s $ with respect to $ p $ are equivalent to the turning points with respect to $ x_2^s $.  General stable manifold
curves will not display such behavior (and indeed, neither does this, if taking more negative $ p $ values or increasing
$ \eps $ further).  We propose as a more general way of determining the folding points the {\em points at which the
curvature exhibits a marked maximum}.  We illustrate the usage of this criterion in Fig.~\ref{fig:curvature}, computed
at the same parameter values as Fig.~\ref{fig:foldpointsp}, in which we show the logarithm of the curvature in
terms of $ p $ (using (\ref{eq:curvature})) and also arclength (using also (\ref{eq:arclength})).
 The same four points
identified in Figs.~\ref{fig:foldpointsp} and \ref{fig:foldpoints} are shown in this figure.  When proceeding from right to
left, i.e., from the hyperbolic trajectory near $ (1,0) $ in Fig.~\ref{fig:foldpoints}, we can see that the curvature is initially
close to zero, corresponding to the almost straight line emanating from the hyperbolic trajectory. Then, the first 
[blue] foldpoint emerges as a local maximum point in curvature.  The next high-curvature points have increasingly 
larger values, and also increasingly sharper peaks, in the curvature plot of Fig.~\ref{fig:curvature}.

Note that in Fig.~\ref{fig:curvature} we used different $ p $ range than in Figs.~\ref{fig:foldpointsp} and ~\ref{fig:foldpoints} to capture more peak points of curvature. We computed the arclength in Fig.~\ref{fig:curvature} with respect to the point $p=41.7151$ and the absolute value of arclength is used for the x-axis. According to the Fig.~\ref{fig:curvature}(a), we can observe more extremes of curvature as $p$ move towards to negative infinity. From the Fig.~\ref{fig:curvature}(b), we can notice that as we proceed along the arclength, the peaks of curvature become larger and the space between two nearest peak points is increased.

\begin{figure}[tb]
	\centering
	\subfigure[]{\includegraphics[width=0.49 \textwidth, height=0.25 \textheight]{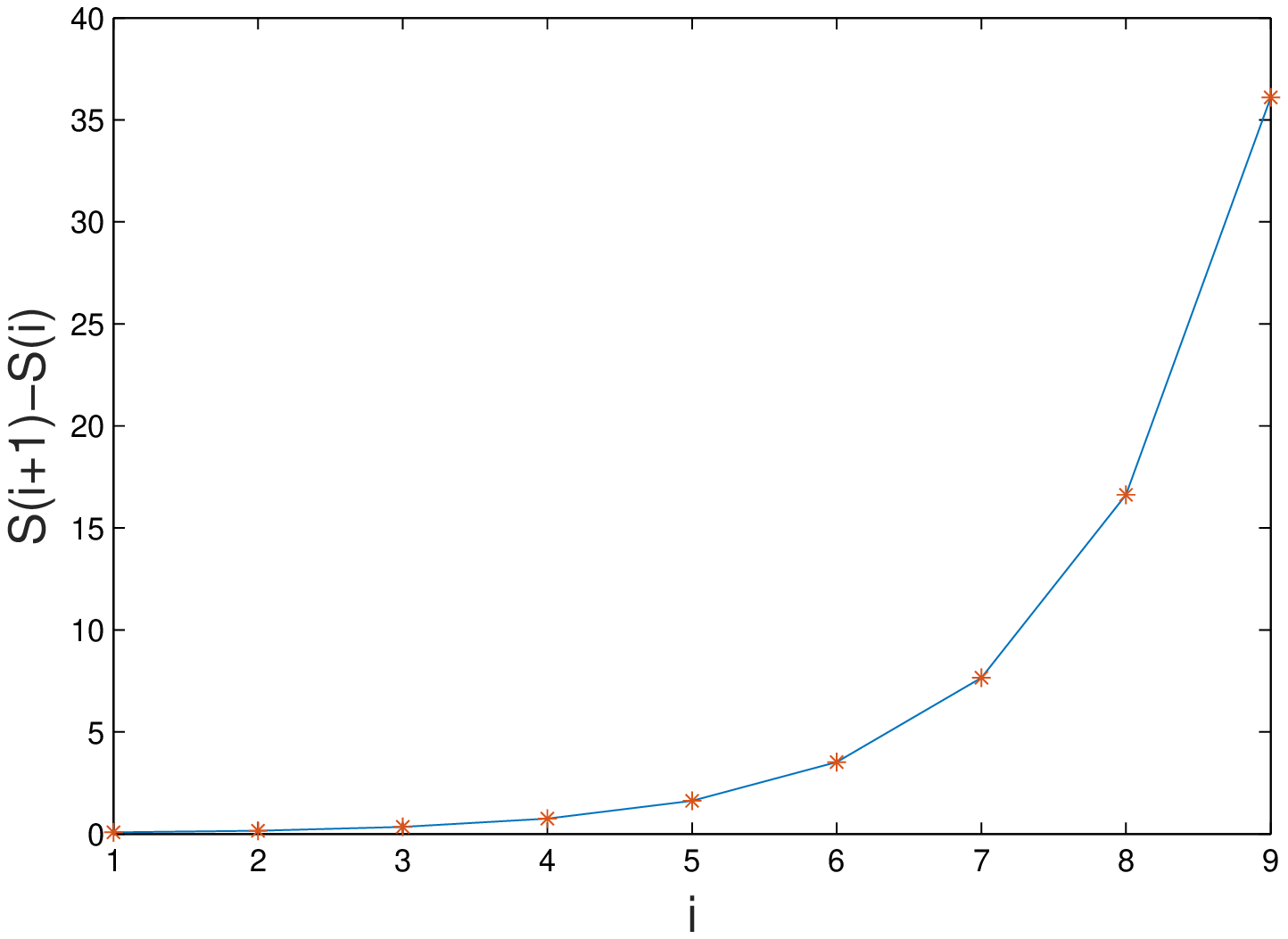}} 
	\subfigure[]{\includegraphics[width=0.49 \textwidth, height=0.25 \textheight]{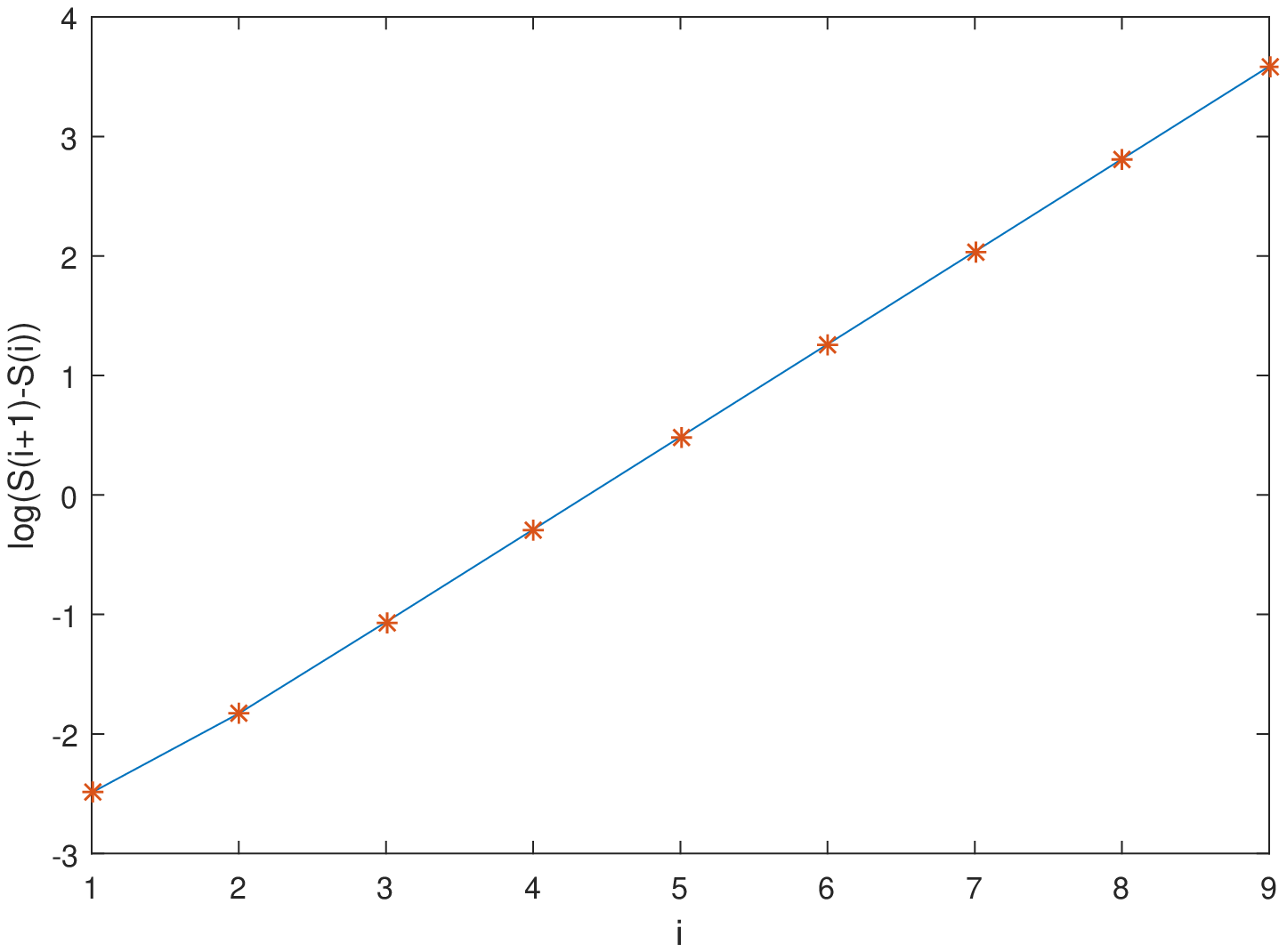}}
	\caption{ (a) The difference of arclength at two consecutive folding points plotted against the $i^{\textnormal{th}}$ folding point, and (b) The logarithm of the difference of arclength at two consecutive folding points plotted against the $i^{\textnormal{th}}$ folding point, at the same parameter values as in Figs.~\ref{fig:foldpointsp} , \ref{fig:foldpoints} and \ref{fig:curvature}.
		}
	\label{fig:difference_arc_length}
\end{figure}

We used ten consecutive folding points (the first at $p=-0.2417$, and the tenth at $p=-0.9485$) from the range $(-1,0)$ to create Fig.~\ref{fig:difference_arc_length}. In Fig.~\ref{fig:difference_arc_length}, the arclength$(S_{i})$ at the  $i^{\textnormal{th}}$ folding point is computed by taking the integral from the point $p=41.7151$ to the $i^{\textnormal{th}}$ folding point, and here we used absolute value of the arclength. These figures give an idea about the arclength distance between folding points. 
Fig.~\ref{fig:difference_arc_length}(b) fits a perfect line in logarithmic scale with the slope of $m=0.7663$ by using linear regression. We can conclude from these results that the arclength between two consecutive folding points grow exponentially, when $p$ moves toward the negative infinity.

\begin{figure}[t]
\subfigure[$(x_1,x_2) $ full-space]{\includegraphics[width=0.5 \textwidth, height=0.25 \textheight]{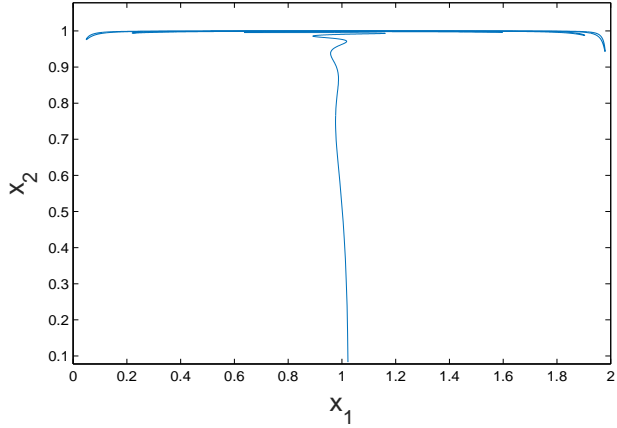}}
\subfigure[Close-up]{\includegraphics[width=0.5 \textwidth, height=0.25 \textheight]{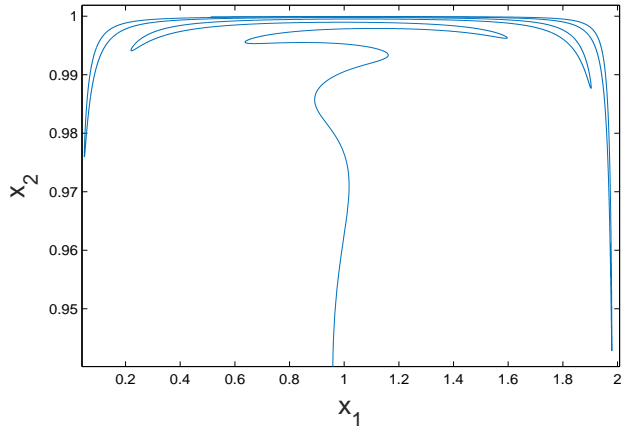}}
 \caption{The stable manifold with the parameter values of the numerical method, N=100000 and $\omega=40$, in (a) $ (x_1,x_2) $-phase space,
 and (b) Zoomed around $x_2 =1 $.}
 \label{fig:stablezoom}
 \end{figure}

\begin{figure}[t]
	\centering
	\includegraphics[width=0.8 \textwidth]{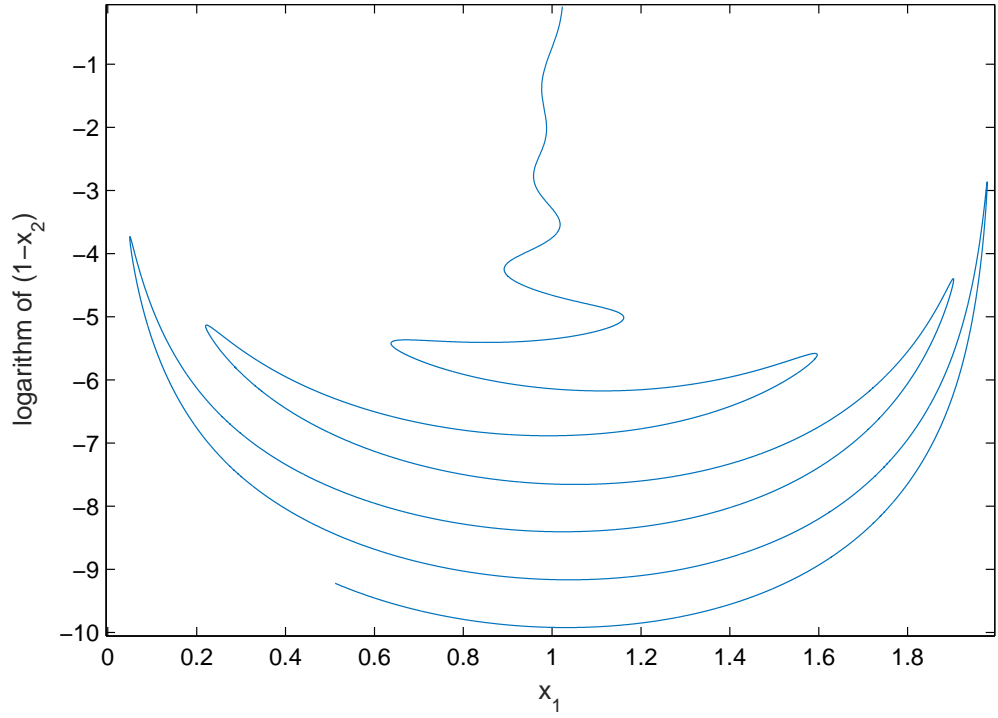}
	\caption{The stable manifold with the parameter values of the numerical method, N=100000 and $\omega=40$, in $ \left( x_1,\ln(1 - x_2) \right) $-space (to elucidate the structure).}
	\label{fig:inverted}
\end{figure}
 
 \begin{figure}[H]
 	\centering
 	\includegraphics[width=0.8 \textwidth]{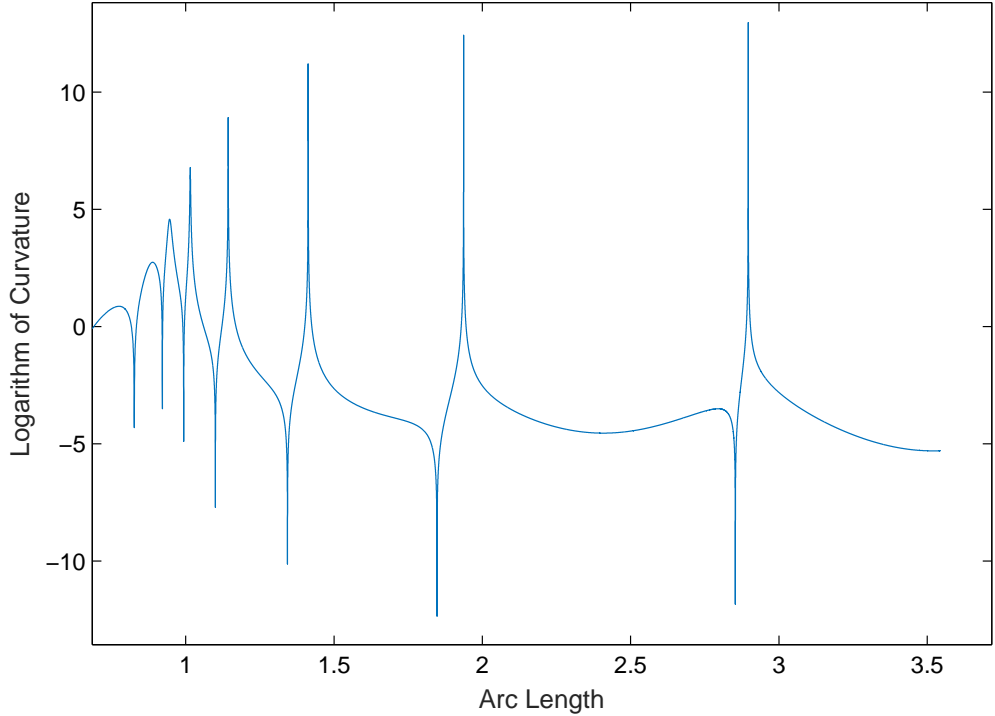}
 	\caption{ \textit{The  logarithm of curvature with respect to the arclength, when $N=100000, \omega=40.$ }}
 	\label{fig7}
 \end{figure} 
 
  From the Fig.~\ref{fig7}, we can notice that as we proceed along arclength, the extremes of curvature become progressively larger due to folded lobes squeezing between lobes as seen in Fig.~\ref{fig:stablezoom}, but spaced progressively further apart, between relatively flat segments, and the near zero curvature points are inflection points.
  
 We note that Figs.~\ref{fig:stablezoom}, \ref{fig:inverted} and \ref{fig7} were plotted numerically and we used same parameter values as in the Figs.~\ref{fig:foldpointsp}, ~\ref{fig:foldpoints},~\ref{fig:curvature} and ~\ref{fig:difference_arc_length} with $N=100000$. Here $N$ is the number of sampling points of the function.
 
From this study, we can conclude that there were infinitely many folding points around $x_{2}=1$ along the stable manifold and there were infinitely many folding points around $x_{2}=0$ along the unstable manifold. We were able to see that the curvature became high at each folding point and these high curvature values were increased when the stable manifold approaches to $x_{2}=1.$ Also we were able to figured out that the arclength between two nearest folding points was increased, when the stable manifold approaches to $x_{2}=1.$

\section{Concluding remarks}

In this paper, we have specifically addressed the concept of folding, particularly in relation to the Double-Gyre flow.  The
folds in the stable and unstable manifolds were used to construct a horseshoe map in this flow, and thereby prove the
implicitly accepted fact that the system is chaotic.  We also show how tracking the curvature is an excellent method for
characterizing folding.  In highlighting the role of folding, we have addressed an aspect of chaos which
is seldom quantified.

\vspace*{0.2cm}
\noindent {\bf {\em Acknowledgments:}} 
SB acknowledges support from the Australian Research Council through grant FT130100484,
and travel support from the University of Adelaide and Clarkson University during his visit to Clarkson when this work
was initiated. EB is supported by, the Army Research Office, N68164-EG and W911NF-12-1- 0276, and the ONR N00014-15-1-2093, and the NGA, and  KGDSP is supported directly by Clarkson University.

\appendix

\section{Proof of Theorem~\ref{theorem:intersect} (Heteroclinic intersections)}
\label{sec:intersect}

The standard method in these situations is to use the Melnikov technique \cite{melnikov,guckenheimerholmes,open},
which is a perturbative technique with respect to $ \eps $.  We will work with the nonautonomous approach as
detailed in \cite{siam_book} in particular, whose approach is to imagine the manifolds as being parametrized by
time $ t $ continuously.  This requires writing (\ref{eq:doublegyre}) in the 
perturbed form
\[
\dot{\vec{x}} = \vec{f} \left( \vec{x} \right) + \eps \vec{g} \left( \vec{x}, t \right) + {\mathcal O}(\eps^2)
\]
in which $ \vec{x} = \left( x, y \right) $, and the  $ {\mathcal O}(\eps^2) $ term is uniformly bounded on $ \Omega \times \R $.
Taylor expansions of (\ref{eq:doublegyre}) enable the identifications
\[
\vec{f}(x,y) = \left( \begin{array}{c}
- \pi A \sin \left( \pi x \right) \cos \left( \pi y \right) \\
\pi A \cos \left( \pi x \right) \sin \left( \pi y \right)
\end{array} \right)
\]
and
\[
\vec{g}(x,y,t) = \left( \begin{array}{ll}
- \pi^2 A \left( x^2 - 2 x \right) \cos \left( \pi y \right) \cos \left( \pi x \right) \\
\pi A \sin \left( \pi y \right) \left[ 2 \cos \left( \pi x \right) \left( x - 1 \right)
- \pi \left( x^2 - 2 x \right) \sin \left( \pi x \right) \right] 
\end{array} \right) \sin \left( \omega t \right) \, .
\]
Here, we follow the methodology of quantifying the signed distance between the perturbed manifolds at a location
$ \bar{\vec{x}}(p) = \left( 1, \bar{x}_2(p) \right) $ on the unperturbed heteroclinic.  At a general time $ t $, the vector from the stable to the
unstable manifold, going through this point and pointing in the $ x_1 $-direction, can be expressed by
\begin{equation}
d(p,t) = \eps \frac{M(p,t)}{\left| \vec{f} \left( \bar{\vec{x}}(p) \right) \right|} + {\mathcal O}(\eps^2) \, , 
\label{eq:distance}
\end{equation}
valid for $ p \in [-P,P] $ and $ t \in [-T,T] $, for $ P $, $ T $ finite.  This is a signed distance, which is positive if the vector
points in the $ + x_1 $-direction, and the Melnikov function in this interpretation is given by \cite{siam_book,open,aperiodic}
\[
M(p,t) = \int_{-\infty}^\infty \vec{f} \left( \bar{\vec{x}}(\tau) \right) \wedge \vec{g} \left( \bar{\vec{x}}(\tau),\tau+t-p \right) \, \d \tau
\]
with the wedge product defined by $ \vec{f} \wedge \vec{g} := f_1 g_2 - f_2 g_1 $ in component form.   Now, substituting the relevant $ \vec{f} $ and $ \vec{g} $, inserting 
$ \bar{\vec{x}}(\tau) = \left( 1, \bar{x}_2(\tau) \right) $, and simplifying leads to
\begin{eqnarray}
M(p,t) & = & \pi^3 A^2  
\int_{-\infty}^\infty  \tanh \left( \pi^2 A \tau \right) \,
\sech \left( \pi^2 A \tau \right) \,\sin \left[ \omega (\tau \! + \! t \! - \! p) \right] \, \d \tau  \nonumber \\
& = & \omega \sech \frac{\omega}{2 \pi A} \sin \left[ \omega (t-p) \right] =: R(\omega) \sin \left[ \omega (t-p) \right] \, , 
\label{eq:melnikov}
\end{eqnarray}
where the final simplification (\ref{eq:melnikov}) is obtainable by writing $ \sin \left[ \omega  \left( \tau + t - p \right)  \right] =
\sin \left[ \omega (t-p) \right] \cos \left( \omega \tau \right) + \cos  \left[ \omega (t-p) \right] \sin \left( \omega \tau \right) $
and splitting the integrals \cite[e.g.]{open,periodic,optimal}, and performing the one non-zero integral that results.  At each
fixed $ t $, $ M(p,t) $ clearly has nonsimple zeros when $ p = t - m \pi / \omega $, $ m \in \Z $.  From (\ref{eq:distance}), this
indicates that $ d(p,t) $ has nearby zeros for small enough $ \left| \eps \right| $.  Thus, the stable and unstable manifolds
intersect---infinitely many times, in fact---in each time slice.  This leads to a heteroclinic tangle near $ x_1 = 1 $, resulting in
complicated transport between the gyres.  This transport can be quantified to leading-order in $ \eps $ 
as an {\em instantaneous flux} \cite{siam_book,open,aperiodic} from the left to the right gyre by 
\begin{equation}
\eps M(p,t) =
\eps R(\omega) \sin \left[ \omega (t-p) \right],
\end{equation}
 as an {\em average flux} \cite{romkedarpoje} by $ \eps R(\omega) $, or
(in terms of lobes created through the intersections) as a {\em lobe area} \cite{romkedar,wiggins} of $ \eps R(\omega) 2 \pi / \omega $.
\hfill $ \Box $

\section{Proof of Theorem~\ref{theorem:stable} (Stable manifold)}
\label{sec:stable}

This proof relies on the expressions for the perturbed stable manifold obtained in Theorem~2.7 in \cite{tangential}
and Theorem~2 in \cite{open}.  While those results are derived for the more general situation of  
{\em compressible} flows which are not necessarily time-periodic, these
relaxations are not necessary in the present context.  Using the notation already introduced in the proof of
Theorem~\ref{theorem:intersect}, the stable manifold at a time $ t $ can be represented parametrically by
\[
\vec{x}^s(p,t) = \bar{\vec{x}}(p) + \eps \left[ \frac{M^s(p,t)}{\left| \vec{f} \left( \bar{\vec{x}}(p) \right) \right|^2} 
\vec{f}^\perp \left( \bar{\vec{x}} \right) + \frac{B^s(p,t)}{\left| \vec{f} \left( \bar{\vec{x}}(p) \right) \right|^2} 
\vec{f} \left( \bar{\vec{x}} \right) \right] + {\mathcal O}(\eps^2)
\]
for $ p \in [P,\infty) $ (any $ P $ which is finite), in which the $ \perp $ notation indicates the rotation of a vector by $ + \pi/2 $,
and expressions for $ M^s $ and $ B^s $ will be given shortly.  In this case, since $ \bar{\vec{x}}(p) $, the unperturbed
stable manifold, simply lies along the line $ x_1 = 1 $, the velocity $ \vec{f} $ along it is directly in the negative $ y $-direction.
Thus, $ \vec{f}^\perp $ is simply the component of $ \vec{f} $ in the $ + x_1 $-direction.  Now, the Melnikov function $ M^s $
is given by
\[
M^s(p,t) = \int_{p}^\infty \vec{f} \left( \bar{\vec{x}}(\tau) \right) \wedge \vec{g} \left( \bar{\vec{x}}(\tau),\tau+t-p \right) \, \d \tau \quad ;
\]
this therefore represents the perturbation of the stable manifold in the {\em normal} direction to the original manifold.  The
tangential perturbation is given by the function
\[
B^s(p,t) := \left| \vec{f} \left( \bar{\vec{x}}(p) \right) \right|^2 \int_0^p  \frac{R^s(\tau) M^s(p, \tau\! + \! t \! - \! p) 
+ \vec{f} \left( \bar{\vec{x}}(\tau) \right) \cdot
\vec{g} \left( \bar{\vec{x}}(\tau), \tau \! +\!  t \! - \! p \right)}{\left| \vec{f} \left( \bar{\vec{x}}(\tau) \right) \right|^2} \, \d \tau \, ,
\]
where
\[
R^s(\xi) := \frac{\left( \vec{f}^\perp \right)^\top \! \left( \bar{\vec{x}}(\xi) \right) 
\left[ \left( D \vec{f} \right)^\top \left( \bar{\vec{x}}(\xi) \right) + \left( D \vec{f}  \right) \left( \bar{\vec{x}}(\xi) \right)
\right] \vec{f} \left( \bar{\vec{x}}(\xi) \right) }{
\left| \vec{f} \left( \bar{\vec{x}}(\xi) \right) \right|^2} \, .
\]
We immediately dispense with the more complicated tangential displacement since it is easy to verify that for the double
gyre, $ R^s \equiv 0 $ and $ \vec{f} \cdot \vec{g} \equiv 0 $.  Therefore, there is no change to the $ x_2 $-coordinate, and
we can write using (\ref{eq:heteroclinic}) that $  x_2^s(p,t) = \bar{x}_2(p) = \frac{2}{\pi} \, \cot^{-1} e^{\pi^2 A p} $.
Using the results derived in the proof of Theorem~\ref{theorem:intersect},
we can write the Melnikov function as
\[
M^s(p,t) = \pi^3 A^2  
\int_{p}^\infty  \tanh \left( \pi^2 A \tau \right) \,
\sech \left( \pi^2 A \tau \right) \,\sin \left[ \omega (\tau \! + \! t \! - \! p) \right] \, \d \tau
\]
which cannot be evaluated in terms of simple functions\footnote{It can be represented in a complicated way in terms of
hypergeometric functions, but this is not particularly illuminating and hence will be avoided.}, unlike in (\ref{eq:melnikov}).
Next, using $ \bar{\vec{x}}(p) $ as given in (\ref{eq:heteroclinic}), we have
\[
\vec{f} \left(  \bar{\vec{x}}(p) \right) = \left( \begin{array}{c} 0 \\
\pi A \sin \left( \pi \bar{x}_2(p) \right) \end{array} \right) = \left( \begin{array}{c} 0 \\
\pi A \sech \left( \pi^2 A p \right) \end{array} \right) 
\]
 and the stable manifold
expression (\ref{eq:stable}) results.  The location of the hyperbolic trajectory can be obtained by appealing directly to
Theorem~2.10 in \cite{tangential}, but here we adopt a more intuitive, formal, approach.  We now have the $ x_1 $-coordinate of
the perturbed stable manifold given by (\ref{eq:stable}), which upon changing the variable of integration (and with the
higher-order term neglected for convenience) can be written as
\[
x_1^s(p,t)  \! = 1 \! + \! \eps \pi^2 A
\! \int_t^\infty \!  \frac{\tanh \left( \pi^2 A (\tau \! - \! t \! + \! p) \right)
\! \sech \left( \pi^2 A (\tau \! - \! t \! + \! p) \right)}{\sech \left( \pi^2 A p \right)}
  \sin \omega \tau \d \tau \, .
 \]
Since the hyperbolic trajectory is approached in the limit $ p \rightarrow \infty $, applying
this limit inside the integral gives
\[
x_1^s(\infty,t) = 1 + \eps \pi^2 A e^{\pi^2 A t} \int_t^\infty e^{-\pi^2 A \tau} \sin \omega \tau \d \tau
\]
which can be integrated and reorganized to give (\ref{eq:xshyp}). The $ x_2 $-coodinate of the hyperbolic trajectory
remains fixed at $ x_2 = 0 $ since it is easy to see that this line is invariant for the full flow (\ref{eq:doublegyre}).  
Next, the reciprocal slope of the manifold at the hyperbolic trajectory is needed.  This is zero when $ \eps = 0 $,
and Theorem~2.2 can be used to prove that the $ {\mathcal O}(\eps) $-correction to this is zero.  More intuitively,
this occurs because
\begin{eqnarray*}
\frac{d x_1^s}{d x_2^s} \Big|_{{\mathbf x}_h^s(t)} & = & \lim_{p \rightarrow \infty} \frac{\partial x_1^s}{\partial p} \Bigg/ \frac{\partial
x_2^s}{\partial p}  \\
& = & \lim_{p \rightarrow \infty}  \frac{\eps \pi^2 A \int_t^\infty \frac{d}{dp} \left[  \frac{\tanh \left( \pi^2 A (\tau \! - \! t \! + \! p) \right)
\! \sech \left( \pi^2 A (\tau \! - \! t \! + \! p) \right)}{\sech \left( \pi^2 A p \right)} \right] \sin \omega \tau \d \tau}{\pi A 
\frac{d}{dp} \sech (\pi^2 A p)} \\
& = & \eps \pi \int_t^\infty \lim_{p \rightarrow \infty} \frac{\frac{d}{dp} \left[  \frac{\tanh \left( \pi^2 A (\tau \! - \! t \! + \! p) \right)
\! \sech \left( \pi^2 A (\tau \! - \! t \! + \! p) \right)}{\sech \left( \pi^2 A p \right) } -1 \right]}{\frac{d}{dp} \sech (\pi^2 A p)}
\sin \omega \tau \d \tau \\
& = &  \eps \pi \int_t^\infty \lim_{p \rightarrow \infty} \frac{\left[  \frac{\tanh \left( \pi^2 A (\tau \! - \! t \! + \! p) \right)
\! \sech \left( \pi^2 A (\tau \! - \! t \! + \! p) \right)}{\sech \left( \pi^2 A p \right) } -1 \right]}{\sech (\pi^2 A p)}
\sin \omega \tau \d \tau \\
& = & \eps \pi \int_t^\infty \lim_{p \rightarrow \infty} \frac{\tanh  \left( \pi^2 A (\tau \! - \! t \! + \! p) \right) - 1}{\sech \left(
\pi^2 A p \right)} \sin \omega \tau \d \tau \\
& = & \eps \pi \int_t^\infty \lim_{p \rightarrow \infty} \frac{ \sech^2 \left( \pi^2 A (\tau \! - \! t \! + \! p) \right) }{- \sech \left(
\pi^2 A p \right) \tanh \left( \pi^2 A p \right)} \sin \omega \tau \d \tau = 0 \, , 
\end{eqnarray*}
where we have utilized the fact that $ \sech \left[ \pi^2 A (\tau - t + p)\right] / \sech \left[ \pi^2 A p \right] \rightarrow 1 $
as $ p \rightarrow \infty $, and L'H\^{o}pital's rule has been used several times.  Thus, the $ {\mathcal O}(\eps) $ 
correction to the slope is zero.  Given that the functions here are all uniformly bounded in suitably high norms, uniformly
for $ t \in \R $, it is clear that the $ {\mathcal O}(\eps^2) $ correction is bounded.
\hfill $ \Box $

\section{Proof of Theorem~\ref{theorem:reentrench} (Fold Re-Entrenchment)}
\label{sec:reentrench}

The outer boundaries of $ \Omega $ can be easily seen to be invariant for any $ \eps $.  For
convenience, we only address the wrapping around ensuing from the right gyre; the left gyre also causes the identical behavior. 
Consider the line $ x_2 = 0 $,
along the interface of the right gyre, that is, for $ x_h^s(t) < x_1 < 2 $.  The flow on this satisfies
\[
\dot{x}_1 = - \pi A \sin \left( \pi \phi(x_1,t) \right) \, . 
\]
Now, when $ \eps = 0 $, we have $ x_h^s(t) = 1 $, and $ \phi $ goes from $ 1 $ at this value to $ 2 $ at $ 2 $.  Since
$ \sin \pi \phi $ is negative in this range, $ \dot{x}_1 $
is positive in this interval.  If $ 0 < \left| \eps \right| < 1/2 $, $ \phi(x_1,t) = 1 $ when
\[
x_1 = \tilde{x} := \frac{ 2 \eps \sin \omega t - 1 + \sqrt{ 1 + 4 \eps^2 \sin^2 \omega t}}{2 \eps \sin \omega t}
\]
as long as $ \sin \omega t \ne 0 $.  (If $ \sin \omega t = 0 $, then $ \tilde{x} = 1 $.)  Thus for $ x_1 \in (\tilde{x},2) $, $ \phi(x_1,t) $ lies between
$ 1 $ and $ 2 $, and therefore the vector field points to the right along the lower boundary of the right gyre, in an interval near
$ x_1 = 2 $.  Using this, and a similar idea for the top of the right gyre, we can obtain the behavior 
as shown in Fig.~\ref{fig:reentrench} by
the blue curves.  This picture is drawn at a general time $ t $, and the red and green represent respectively the unstable and the
stable manifold, whose behavior of this form is guaranteed by Theorems~\ref{theorem:stable} and \ref{theorem:unstable}.  Only
parts of the stable manifold near $ x_2 = 0 $ and the unstable manifold near $ x_2 = 0 $ and $ x_2 = 1 $ is shown. While the
arrows drawn on the blue bounding lines are the instantaneous directions of the velocity, those drawn on the stable/unstable
manifolds are not necessarily the instantaneous velocity directions, since these manifolds, and their anchor points $ \vec{x}_h^{s,u}(t) $, are themselves moving (mostly horizontally in the regions near $ x_2 = 0 $ and $ x_2 = 1 $).  The true instantaneous velocity of
particles on these manifolds is the superposition of the indicated arrows on the manifolds, and this additional motion.

Now, it must be borne
in mind that the unstable manifold intersects the stable one infinitely often near $ (x_h^s(t),0) $, with the intersection points
accumulating towards this instantaneous hyperbolic trajectory location.  However, the lobe structures created as a result of 
this intersection must have equal areas, since under iteration of the Poincar\'{e} map $ P $ which samples the flow from this time
$ t $ to the time $ t + 2 \pi / \omega $ (i.e., strobing the flow at the period of the velocity field), these lobes must get mapped
to one another.  The lobe with end marked by $ B $, must get mapped to the next lobe with end marked by $ P(B) $.  This lobe
must get thin in the $ x_2 $-direction (indeed, this width is almost not discernible in Figure~\ref{fig:reentrench}) because the
intersection points accumulate to $ (x_h^s(t),0) $.  However, the flow of (\ref{eq:doublegyre}) is incompressible, and thus
area-preserving.  The lobe which has $ P(B) $ marked at its end must therefore have the same area as that marked with a
$ B $, and this is only achievable if it extends outwards.  This extension in the $ x_1 $-direction of the unstable manifold is also
implied by the formul\ae{} for $ x_1^u(p,t) $ shown in Theorem~\ref{theorem:unstable}.  One can therefore determine parts of the
unstable manifold which are {\em arbitrarily} close to the line $ x_2 = 0 $, and there will be regions of this manifold which have
$ x_1 $-coordinates greater than $ \tilde{x} $.  By continuity, the velocity at such a location can be made arbitrarily close to the velocity
on $ x_2 = 0 $.  Thus, if considering the blob marked $ B $ in Figure~\ref{fig:reentrench} which is at an end of a lobe structure
(where the unstable manifold folds)
and assuming that this has been chosen to be within this region of influence, as time passes it will get pulled along by a velocity
which is very close to that of the blue lines.  Eventually, therefore, it must get pulled all the way around the right gyre, and
emerge along the blue line at the top of the right gyre.  

While the flow along this boundary line is to the left for $ x_1 $ near $ 2 $,
we want to show something more specific: that the flow along this line approaches the hyperbolic trajectory location $ x_h^u(t) $,
as approximated in (\ref{eq:xuhyp}).  Focus, then, on flow along this blue line, that is on the invariant line $ \left\{ x_2 = 1, 0 < x_1 < 2 
\right\} $, which obeys
\[
\dot{x}_1 = \pi A \sin \left( \pi \phi(x_1,t) \right) \, . 
\]
Let $ z(t) = x_1(t) - x_h^u(t) $, and suppose that $ x_1(0) > x_h^u(0) $.  Since trajectories cannot cross on this one-dimensional
phase space, it is clear that $ x_1(t) > x_h^u(t) $ for $ t > 0 $.  Now $ 0 < x_h(t) < x_1(t) < 2 $  because 
the end points $ 0 $ and $ 2 $ are fixed points of the above, even if the flow is nonautonomous.  Therefore $ 0 < z(t) < 2 $, 
and $ z(t) $  can be 
shown to satisfy the differential equation
\begin{eqnarray*}
\dot{z} & = & 2 \pi A \sin \left[ \frac{\pi z}{2} + \frac{\eps \pi z \sin \omega t}{2} \left( x_1 + x_h^u -2 \right) \right] \\
& & \times \cos \left[  \frac{\pi (x_1 + x_h^u)}{2} + \frac{\eps \pi (x_1 +x_h^u) \sin \omega t}{2} \left( x_1 + x_h^u - 2 \right) \right] \\
& = & - 2 \pi A \sin \left[ \frac{\pi z}{2} + \frac{\eps \pi z \sin \omega t}{2} \left( x_1 + x_h^u -2 \right) \right] \\
& & \times \sin \left[  \frac{\pi (x_1 + x_h^u-1)}{2} + \frac{\eps \pi (x_1 +x_h^u) \sin \omega t}{2} \left( x_1 + x_h^u - 2 \right) \right] 
\, .
\end{eqnarray*}
When $ \eps = 0 $, $ x_h^u(t) = 1 $, and in this situation 
\[
\dot{z} = - 2 \pi A \sin \frac{\pi z}{2} \sin \frac{\pi x_1(t)}{2} < 0 \, , 
\]
whose velocity field is sign definite since both $ x_1 $ and $ z $ must lie in $ (0,2) $.  Despite being nonautonomous, its solution 
$ z $ must decay to
the fixed point $ z = 0 $.  When $ \eps \ne 0 $, noting also that $ x_h^u = 1 + {\mathcal O}(\eps) $, it is clear that one can find
$ \left| \eps \right| $ small enough such that the sign definite nature will persist if $ x_1(t) $ were chosen sufficiently close to $ 
x_h^u(t) $.  Therefore, along the blue line at the top of Fig.~\ref{fig:reentrench}, for suitably small $ \left| \eps \right| $, trajectories
will be attracted towards the hyperbolic trajectory $ x_h^u(t) $.  

Once we have this property, continuity ensures that trajectories inside $ \Omega $ but near to this must also follow the behavior of
proceeding towards the left.  Trajectories can be made to approach $ x_h^u(t) $ arbitrarily closely, by choosing trajectories which
were sufficiently close to the blue line.  However, the fluid blob $ B $ as shown in Fig.~\ref{fig:reentrench} will at some time
in the future be as close as we like to the `heteroclinic network' shown in blue in Fig.~\ref{fig:reentrench}, and thus will eventually
be subject to behavior imputed for the line $ x_2 = 1 $.  When approaching $ x_h^u(t) $, this blob will therefore be subject to the
unstable manifold emanating from $ x_h^u(t) $, and get pulled down along it.  

Consider the point $ u $, which is at the leading-edge of the lobe marked $ B $ in Figure~\ref{fig:reentrench}. That is, this is a fold
point.  By the above argument, the flow of the Double-Gyre will ensure that its image $ u' $ will eventually be within $ N_\delta $.  
We want to show the existence of a fold point within $ N_\delta $.  However, there is no guarantee that 
the point $ u' $ will also correspond to a leading-edge, i.e., a fold point.  To establish the existence of a fold point, we argue that
the unstable manifold must pass through $ u' $.  Now, both ends of the unstable manfold must wrap back all the way around the
boundary of $ \Omega $, adjacent to the blue lines, and come back to intersect the stable manifold of $ \left( x_h^s(t), 0 \right) $
near to this point, since these intersection points accumulate towards $ \left( x_h^s(t), 0 \right) $.   Hence, both ends of the
unstable manifold which pass through $ u' $ must come all the way back.  This ensures that there must be a fold point within
$ N_\delta $; the unstable manifold must `bend back' to achieve this.

It is also instructive to think of what happens in terms of the Poincar\'{e} map and lobes.  It has been argued that there are 
infinitely many lobes `below' the one pictured near $ B $ in Figure~\ref{fig:reentrench}.  As one proceeds `downwards,' each
of these lobes is closer to the blue heteroclinic network than the previous one, and therefore subject to the motion along the
network more.  Thus, each successive lobe will get elongated along the network more.  The end result from this process is
that in $ N_\delta $, in addition to the lobe $ L $ pictured in Figure~\ref{fig:reentrench}, there will be an infinite number of lobes which 
accumulate towards the unstable manifold emanating from $ x_h^u $.  These will stretch along the unstable manifold (they cannot
intersect because the lobe boundaries are themselves part of the {\em same} unstable manifold), and therefore will follow the
undulations that the primary part of the unstable manifold has been shown to have.
\hfill $ \Box $

\section{Expressions from the analytical approximations}
\label{sec:approx}

Here we list some expressions related to determining the curvature and fold points from the analytical
approximations given by Theorems~\ref{theorem:stable} and \ref{theorem:unstable}.  At each fixed time
$ t $, the primary stable/unstable manifold curves can be thought of as being given parametrically by
(\ref{eq:stable}) and (\ref{eq:unstable}), where $ p $ is the parameter.  We will only show calculations for the
stable manifold, since the unstable manifold calculations are similar.  By applying integration by parts
and a straightforward change-of-variable to (\ref{eq:stable}), the $ {\mathcal O}(\eps) $ expression
for the $ x_1 $-coordinate of the stable manifold, given in (\ref{eq:stable}), 
can be recast as
\[
x_{1}^{s}(t) =1 +\eps\sin(\omega t) +\eps \omega \cosh(\pi^2 A p) \int_{t}^{\infty} \sech \left[ \pi^2A(u-t+p)
\right] \cos(\omega u)
\, \d u \, .
\]
Its derivative is therefore
\begin{eqnarray*}
 \frac{dx_1^s}{dp} & = & \eps \omega \pi^2 A \cosh(\pi^2Ap)  \int_{t}^{\infty} \sech \left[ \pi^2A(u-t+p) \right] \\
 && \times   \left\{ \tanh \left[ \pi^2 A p \right] - \tanh \left[ \pi^2A(u-t+p) \right] \right\} \cos(\omega u) \, \d u \, . 
 \end{eqnarray*}
While not obvious in the above representation, it turns out that $ d x_1 / d p $ takes on a sinusoidal form in $ t $, 
which helps us locate its zeros quickly.  To obtain this form, we first define
\[
f_1(v,p) = \eps \omega \pi^2 A \cosh(\pi^2Ap) \sech(\pi^2Av) \left[ \tanh(\pi^2 A p) - \tanh(\pi^2 A v) \right] \, , 
\]
and
\[
J(p)=\sqrt{\left(\int_{p}^{\infty} f_{1}(v,p) \cos(\omega v) \d v\right)^{2}+\left(\int_{p}^{\infty} f_{1}(v,p) \sin(\omega v) \d 
v\right)^{2}} \, .
\]
Then, after some trigonometric manipulations, it is possible to write
\begin{equation}
\frac{d x_1}{d p} = J(p) \cos \left[ \omega (p-t)-\theta(p) \right] \,  , \,  
\theta(p) = \cos^{-1}\left(\frac{ \int_{p}^{\infty} f_{1}(v,p) \cos(\omega v) \d v}{J(p)}\right) \, ,
\label{eq:x1p}
\end{equation}
from which zeros can be obtained easily using a Newton-Raphson method.  These represent parameter values $ p $
corresponding to fold points, as long as $ d^2 x_1 / d p^2 $ is sign definite.  This takes the form
\begin{small}
\begin{eqnarray}
\frac{d^2 x_1^s}{d p^2} & = & \eps \omega \pi^4 A^2 \cosh(\pi^2 A p) \int_{t}^{\infty} \d u \sech \left[ \pi^2A(u-t+p) \right] \cos(\omega u) \nonumber \\
& & \times \left\{  2 \tanh^{2} \left[ \pi^2A(u-t+p) \right] -\tanh \left[ \pi^2Ap \right] \tanh \left[ \pi^2A(u-t+p) \right]
-\tanh^{2} \left[ \pi^2Ap\right]  \right\} \nonumber \\
& & + \eps \omega \pi^4 A^2  \sinh(\pi^2Ap) \int_{t}^{\infty} \d u \sech\left[ \pi^2A(u-t+p)\right] \cos(\omega u) \nonumber \\
& & \times \left\{ \tanh \left[ \pi^2Ap \right] - \tanh \left[ \pi^2A(u-t+p))
\right] \right\} \, .
\label{eq:x1pp}
\end{eqnarray} 
\end{small}
It is straightforward to compute the $ p $-derivatives of the $ x_2 $-coordinate in (\ref{eq:stable}) to be
\begin{equation}
\frac{dx_2^s}{d p}=-\pi A \sech \left( \pi^2 A p \right)  \quad , \quad
\frac{d^2 x_2^s}{d p^2} = \pi^3 A^2 \sech \left( \pi^2 A p \right) \tanh \left( \pi^2 A p \right) \, .
\label{eq:x2pp}
\end{equation}
Given expressions (\ref{eq:x1p}), (\ref{eq:x1pp}) and (\ref{eq:x2pp}), the following geometrical quantities are
easy to compute:
\begin{itemize}
\item The arclength between two points with parametric coordinates $ p_1 $ and $ p_2 $:
\begin{equation}
S(p_1,p_2) = \int_{p_1}^{p_2} \sqrt{\left(\frac{dx_1^s}{d p}\right)^{2}+\left(\frac{dx_2^s}{d p}\right)^{2}} \d p \, .
\label{eq:arclength}
\end{equation}
\item The curvature at a general location $ p $ on the stable manifold:
\begin{equation}
\kappa(p) =\frac{\left| \frac{d^2 x_2^s}{d p^2} \frac{d x_1^s}{d p} - \frac{d^2 x_1^s}{d p^2} \frac{d x_2^s}{d p} \right|}
{ \left[  \left( \frac{d x_1^s}{d p} \right)^2 + \left( \frac{d x_2^s}{d p} \right)^2 \right]^{3/2} } \, .
\label{eq:curvature}
\end{equation}
\end{itemize}
These expressions can also be used to determine the arclength and curvature of the {\em unstable} manifold, by substituting
the expressions for $ x_{1,2}^u $ instead.

\bibliographystyle{plain}
\bibliography{doublegyre1}

\end{document}